\documentclass{article}

\usepackage{arxiv}
\usepackage{amsmath}
\usepackage{graphicx}
\usepackage{epstopdf}
\usepackage[utf8]{inputenc} 
\usepackage[T1]{fontenc}    
\usepackage{hyperref}       
\usepackage{url}            
\usepackage{booktabs}       
\usepackage{amsfonts}       
\usepackage{nicefrac}       
\usepackage{microtype}      
\usepackage{lipsum}
\usepackage{algorithm}
\usepackage{algorithmic}
\usepackage{booktabs}
\usepackage{lscape}
\usepackage{caption}
\usepackage{subcaption}

\newtheorem{defi}{Definition}[section]
\newtheorem{thm}{Theorem}[section]

\title{Analysis and Simulation of a Coupled Diffusion based 
Image Denoising Model}

\author{
 Subit K. Jain \\
  School of Basic Sciences\\ 
  Indian Institute of Technology Mandi\\
  PIN 175005, INDIA\\
  \texttt{jain.subit@gmail.com} \\
    \And
 Sudeb Majee \\
  School of Basic Sciences\\ 
  Indian Institute of Technology Mandi\\
  PIN 175005, INDIA\\
  \texttt{sudebmajee@gmail.com} \\
    \And
  Rajendra K. Ray \\
  School of Basic Sciences\\ 
  Indian Institute of Technology Mandi\\
  PIN 175005, INDIA\\
  \texttt{rajendra@iitmandi.ac.in} \\
    \And
   Ananta K. Majee \\
   Department of Mathematics\\
   Indian Institute of Technology Delhi\\
   PIN 110016, INDIA \\
  \texttt{majee@maths.iitd.ac.in} \\ 
}

\begin{document}
\maketitle

\begin{abstract}
In this study, a new coupled Partial Differential Equation (CPDE) based image denoising model incorporating space-time regularization into non-linear diffusion is proposed. This proposed model is fitted with additive Gaussian noise which performs efficient image smoothing along with the preservation of edges and fine structures. For this purpose, we propose a new functional minimization framework to remove the image noise, which results in solving a system of three partial differential equations (PDEs). Our proposed model is dissimilar from the existing CPDE models as it includes two additional evolution equations to handle edge strength function and data fidelity term. These two evolution equations control the smoothing process and force the resultant denoised solution to be close to the initial solution. To the best of our knowledge, the proposed model is the only work, which deciphers the combined effect of both the terms using separate PDEs. Furthermore, we establish the existence and uniqueness of a weak solution of the proposed system using the time discretization method with $H^1$ initial data. Finally, we used a generalized weighted average finite difference scheme to efficiently solve the coupled system and experiment results show the effectiveness of the proposed CPDE model.
\end{abstract}

\section{Introduction}
Reduction of generated noise in the restored images without loss of spatial resolution is highly desirable in many image processing applications. The behavior of noise reduction methods is entirely different for large areas with well-defined borders as well as for poorly defined borders. It has been observed that the goal of noise reduction with feature preservation and applicability to the multiple acquisition systems is difficult to achieve through a single method \cite{alvarez1993axioms,gonzalez2002digital}. Therefore, the accurate reconstruction of key features in an image from the noisy measurement is the main objective of our study.

The noise caused by image acquisition devices is often modeled by Gaussian random distribution. Generally, the degraded image can be modeled as
\begin{equation}
J = I + \eta,
\end{equation}
where $J$ is the noisy image, degraded by Gaussian noise with zero mean, defined on an image space $\Omega$, and $I$ depicts the clean image, and $\eta$ is the additive noise function. Our main task is to obtain the reconstructed image from the observed noisy image.
 
In recent years, PDE based methods have become widely applicable for noise removal and signal reconstruction, due to their well-established mathematical properties and their general wellposedness
\cite{alvarez1993axioms,aubert2006mathematical,mikula2002image,nolen2009partial,weickert1998anisotropic,weickert2001applications,witkin1984scale}. Obtaining the steady-state solution is the basic objective of PDE based image restoration problem. This goal can be achieved through two approaches. In the first approach, we directly find a steady state solution of the diffusion model. Whereas in the second approach, the steady-state solution is obtained from the evolution equation of the Euler-Lagrange equation associated with the energy functional. Several techniques are employed to handle the challenges posed by restoration problems. The availability of vast literature does not allow us to review all of them. However, we briefly report some of the robust techniques such that anisotropic diffusion \cite{aubert2006mathematical,weickert1998anisotropic,weickert2001applications,catte1992,perona1990scale,weickert1997review}, fourth order diffusion
\cite{liu2011adaptive,siddig2018image,you2000fourth}, complex diffusion \cite{araujo2012stability,araujo2015stability} and variational based methods \cite{chambolle1997image,elliott2009numerical,rudin1992nonlinear,prasath2010hybrid,tsai2005total,zanella2018serial}.
In all of the above-mentioned methods, spatial regularization is used for diffusion function, which is not able to inject the past information into the iterative process of diffusion. Therefore, the scheme may not be effective for solving image denoising problems which contain huge textures in their image domain. To handle this, applications of time-delay regularization with spatial regularization are also available in the literature \cite{amann2007time,belahmidi2005time,guo2011reaction,luo2006coupled,nitzberg1992nonlinear,prasath2014system}. The space-time regularization is more suitable approach than Gaussian filtering as it enables to incorporate the information obtained along with the scales into the diffusion process. This type of regularization is used to recover fine structures, especially for images degraded with higher noise levels.
Hence apart from a single PDE based approaches, several studies are utilizing coupled PDE to perform image restoration \cite{guidotti2015anisotropic}. In this regard, Nitzberg and Shiota \cite{nitzberg1992nonlinear} introduce a delay in time to calculate diffusivity term as,
\begin{equation}
I_t  = \nabla\left(\frac{1}{1+u^2}\nabla I\right),  \hspace{2.6cm} \text{in} \hspace{0.2cm} \Omega \times (0,T),
\end{equation}
where $u$ in diffusion function is calculated using the relation:  $u_t =\varphi(G_{\xi} \ast|\nabla {I}|^{2}-u)$ and $T$ is a fixed time. Here, $\xi$ is a positive smoothing parameter required for the analytical well-posedness of the model. So, the term $u$ updates the edge information from the gradient of the updated image and past information of itself. Hence, instead of spatial regularization, $u$ represents a temporal regularization in the formation of the diffusion term. This coupled PDE system is quite similar to the Perona-Malik model \cite{perona1990scale} as there is no spatial smoothing term. In \cite{nitzberg1992nonlinear}, the authors show that the coupled PDE system admits a unique classical solution $(I,u)$ in any dimension and satisfy the maximum principle. They also explain that the system does not produce spurious features. A similar idea of coupled PDE system using time-delay regularization is proposed by Luo et al. \cite{luo2006coupled}, for image denoising. They propose the estimation of a better edge map by substituting the isotropic diffusion by a nonlinear diffusion equation, to produce a regularized version of the diffusion function for image smoothing. Subsequently, Guo et al. \cite{guo2011reaction} propose and analyze the well-posedness of a reaction-diffusion system for image denoising. They deploy the $H^{-1}$ norm instead of $L1/L2$ norm with the total variation model \cite{osher2003image}, to preserve the oscillatory and texture patterns, as follows,
\begin{equation}
E_{H^{-1}}= \int_{\Omega} |\nabla I|\,d \Omega \ + \lambda ||J-I||_{H^{-1}(\Omega)}^2.
\end{equation}
Here, $||.||_{H^{-1}(\Omega)}^2=  \int_{\Omega} |\nabla \Delta^{-1}(.)|^2\,d\Omega$. This model is efficient for denoising and decomposing the image into cartoon plus texture. Now evolution equation of the Euler-Lagrange equation of the above functional, in which two PDEs are interacting with each other to calculate noise fidelity term, can be given as,
\begin{equation}
\left.\begin{aligned}
I_t &= \nabla\left(\frac{\nabla I} {|\nabla I|}\right)-2 \lambda v,  \hspace{2.1cm} \text{in} \hspace{0.2cm} \Omega \times (0,T),  \\
v_t &= \Delta v-(J-I),  \hspace{2.5cm} \text{in}\hspace{0.2cm} \Omega \times (0,T). 
\end{aligned}\right\}
\end{equation}
In the above model, another PDE is used to calculate the fidelity $\left( v \right)$ term efficiently. The main difficulty with this model is that the term $\frac{1}{|\nabla I|}$ is singular when $|\nabla I|=0$. However, the numerical simulations involve the regularization of $|\nabla I|$ with $|\nabla I|+\varepsilon$. A good survey on the existing CPDE approaches can be found in \cite{guidotti2015anisotropic}.

The objective of the present work is to systematically develop a new framework of coupled nonlinear diffusion model. The present approach differs from previous existing CPDE based works as it utilizes separate evolution equations to calculate edge strength function as well as the data fidelity term. The originality of the proposed coupled non-linear PDEs algorithm lies in the combined effect of:
\begin{itemize}
\item[i.] diffusion coefficient with space-time regularization instead of space regularization only;
\item[ii.]  $H^{-1}$ norm based time-varying fidelity term.
\end{itemize}
Moreover, due to the interest in theoretical research, we establish the well-posedness of the proposed CPDE model \eqref{maina}-\eqref{maine}. To prove the well-posedness, we use a time discretization method along with essential a-priori estimates and classical results of compact inclusion in Sobolev spaces \cite{raadams1975}; by utilizing some important inequalities, lemmas and theorems \cite{lcevans1998} and then tends the time discretization parameter to zero. This imparts a major role of this study as it is very important for the numerical computation.
Further, to simulate the image denoising, we apply a higher order accurate generalized weighted average finite difference scheme with an advanced iterative solver (Hybrid Bi-Conjugate Gradient Stabilized method) to solve the algebraic system of equations generated from the discretization process \cite{jain2015Elsevier}. The image quality of denoised images using the proposed CPDE model is compared against several existing non-linear diffusion models, variational based models and CPDE based models. In practice, our calculations indicate that the proposed approach allows effective image smoothing on fine scales even when the images are degraded with a higher level of noise.

The rest of the paper is arranged in the following sequence. Section \ref{sec:proposed} describes the mathematical formulation of the proposed coupled nonlinear diffusion model. In section \ref{sec:weaksolutions} we establish the existence and uniqueness of a weak solution to the proposed PDE model. Section \ref{sec:numerical} shows an appropriate numerical realization in terms of higher order accurate implicit finite difference scheme. Numerical validation of the proposed approach using experimental study is carried out in section \ref{sec:result}. At the end, we summarize our observations in section \ref{sec:conclusion}.

\section{Proposed Coupled Diffusion Model}
\label{sec:proposed}
The main aim of this work is to present a novel approach for image smoothing by the evolution of coupled non-linear PDEs. Among the existing CPDE based models mentioned in the previous section, the choices of diffusion function and fidelity term include the space and time regularization. Numerical simulations and experimental results depict that the space-time regularization plays a important role in the quality of recovered images degraded by noise. We note that all of these existing coupled PDEs described in the last section utilizes the time-delay regularization to calculate either diffusion coefficient or fidelity term, and sometimes fail to preserve significant contents of the image. Thus, one issue of these smoothing algorithms is how to design a non-linear diffusion algorithm, which can remove the noise and preserve the small features, simultaneously. For this purpose, inspired by the impressive performance of time-delay regularization in image smoothing, we proposed the following energy functional:
\begin{equation}\label{energy}
\text{min}E_u(I)=\int_{\Omega} \frac{1}{2}g(u)|\nabla I|^2\,d \Omega + \lambda||I_0-I||_{H^{-1}(\Omega)}^2, 
\end{equation}
where $I_0$ represents the high pass filtered version of image $J$ and $\lambda$ is weight parameter.
Generally, $\lambda$ is chosen to be inversely proportional to the variance of the noise in the given image. The first term represents a regularizing term producing a smooth varying variable function. The term $g(u)$ serves the purpose of detecting the edges in the image, and $u$ is the edge variable. The second term forces the denoised solution $I$ to be a close approximation to $I_0$ and referred to as data fidelity term. In the data fidelity term, we have replaced $L^2$ norm of ($I_0-I$) by $H^{-1}$ norm to preserve the oscillatory and texture information more appropriately in a denoised image. The solution of the above minimization problem (\ref{energy}) can be given by the steady-state of the Euler-Lagrange equation for $E_u(I)$:
\begin{equation}
 \frac{1}{2 \lambda} \nabla (g(u)\nabla I) =  \Delta^{-1}(I_0-I).
\end{equation}
This is equivalent to following equations,
\begin{equation}\label{rd}
\left.\begin{aligned}
\nabla (g(u)\nabla I)-2\lambda v = 0,\\
\nabla(\nabla v)-(I_0-I) = 0.
\end{aligned}\right\}
\end{equation}
To this end, we consider the function $u$ as the following, 
\begin{equation}\label{ns}
u_t=\alpha(|\nabla {I_\xi}|^{2}-u+\frac{\beta^2}{2}\Delta u),
\end{equation}
where $\alpha >0 $ and $\beta > 0$ are parameters to be specified; average time delay ($\alpha^{-1}$) and amount of the spatial smoothing ($\beta$) control the deformation and smoothness of the edge information and chosen as described in \cite{nitzberg1992nonlinear}. And, $I_{\xi}=G_{\xi}(X) \ast I$ where two dimensional Gaussian function $G_{\xi}(.)$ at each pixel, $X \in \Omega$ is adopted as,
\begin{equation}\label{eq:gaussian_function}
G_{\xi}(X)=\left( \xi \sqrt{2 \pi} \right)^{-1} \text{exp} \left( \frac{-|X|^2}{2 \xi^2} \right),
\end{equation}
where $``\ast "$ denote the convolution operator and $\xi>0$ is the standard deviation of the function. 
Finally, our proposed model in which a set of PDEs interacting with each other, can be expressed as follows,
\begin{align}\label{maina}
\dfrac{\partial I}{\partial t}  &= \nabla(g(u)\nabla I)-2 \lambda v,  \hspace{2.6cm} \text{in} \hspace{0.2cm} \Omega_T:=(0,T)\times \Omega, \\
\label{mainb}
\dfrac{\partial u}{\partial t} &=\alpha(h(|\nabla {I_\xi}|^{2})-u+\frac{\beta^2}{2}\bigtriangleup u),  \hspace{1.2cm} \text{in} \hspace{0.2cm}  \Omega_T, \\
\label{mainc}
\dfrac{\partial v}{\partial t} &= \nabla(\nabla v)-(I_0-I),  \hspace{2.6cm} \text{in}\,\,\, \Omega_T, \\
\label{maine}
I(0,x)&=I_0(x), \hspace{0.5cm} v(0,x)=0, \hspace{0.5cm}	u(0,x)=G_{\xi}\ast|\nabla I_0|^{2}, \hspace{0.5cm} \text{in} \hspace{0.2cm}\Omega.
\end{align}
In the above model, diffusion coefficient $g(u)$ is chosen as
\begin{equation}
\label{eq:gv}
g(u)=\dfrac{1}{1+\frac{|G_{\xi}\ast u|}{k^2}},
\end{equation}
where $k>0$ is a threshold parameter, and $u$ represents edge strength at each iteration. Note that we replace $|\nabla {I_\xi}|^{2}$ in \eqref{ns}  by  $h(|\nabla {I_\xi}|^{2})$, where $h$ being a sort of truncation, see \cite{belahmidi2005time}. Our proposed model is different from the models reported earlier in \cite{guidotti2015anisotropic}, as it yields two separate evolution equations (PDEs) to handle the diffusion coefficient and data fidelity term. These two evolution equations are responsible for stopping the image smoothing at the edges as well as textures and forces the resulting solution to be a close approximation to the given initial image. The proposed model is in a spirit similar to the reaction-diffusion model equation in \cite{guo2011reaction}; however, the diffusion function is very different in our case. In \eqref{maina}-\eqref{maine}, the smoothing equation \eqref{maina} has different edge variable $u$ and fidelity term $v$, obtained from two different PDEs. The data fidelity term between $I$ and $I_0$ can be handled by function $v$, which can be obtained from equation \eqref{mainc}. Whereas the function $u$ in the diffusion coefficient can be calculated from equation \eqref{mainb}. In summary, the proposed model achieves a suitable edge map and fidelity between noisy image and restored image at each iteration, which ultimately leads to quality denoising results. 
Hence, the proposed model provides a potential approach for image noise removal and enhancement.
\section{Existence and uniqueness of weak solution}
\label{sec:weaksolutions}
In this section, we establish the well-posedness of \eqref{maina}-\eqref{maine} with Dirichlet's boundary condition for the fidelity variable and Neumann's boundary condition 
for other two variables. Furthermore, for simplicity we choose all the constants involved in the equations \eqref{maina}-\eqref{maine} equals to $1$. Note that, $h:\mathbb{R}^+ \rightarrow \mathbb{R}^+$ is a
bounded, Lipschitz continuous function with Lipschitz constant $c_h$ such that
\begin{align}\label{bound:h}
 \delta \leq h(\tilde{u})\leq 1\, \quad \forall\, \tilde{u}\in \mathbb{R}^+\,, 
 \end{align}
for some $\delta >0$. Also from \eqref{eq:gv} we observe that, $g: \mathbb{R} \rightarrow \mathbb{R}^+ $ is a bounded, decreasing and Lipschitz continuous function with Lipschitz constant $\frac{C_{\xi}}{k^2}$. Moreover, $g(0)=1$ and $\underset{u \rightarrow + \infty }\lim g(u)=0$.
\subsection{Technical framework}
We denote by $H^k(\Omega)$, $k$ is a positive integer, the set of all functions $I:\Omega \rightarrow \mathbb{R}$ such that $I$ and its distributional derivatives $\frac{\partial^m I}{\partial x^m}$
of order $|m|=\sum_{j=1}^{2}m_j \leq k$ all belongs to $L^2(\Omega)$. $H^k(\Omega)$ is a Hilbert space endowed with the norm
\begin{align*}
 \|I\|_{H^k}:=||I||_{H^k(\Omega)}=\Big(\underset{|m|\leq k }{\sum}\int_{\Omega} \Big|\dfrac{\partial^m I}{\partial x^m}\Big|^2\,dx \Big)^{1/2}\,.
 \end{align*}
We denote by $L^{p}(0,T;{H^k(\Omega)})$, $p>1$ and $k$ is a positive integer, the set of all measurable functions $I:[0,T] \rightarrow H^k(\Omega)$ such that 
$$\|I\|_{L^p(0,T;H^k(\Omega))}:=\Big(\int_{0}^{T} || I(t)||^p_{{H^k(\Omega)}}\,dt \Big)^{1/p} < \infty.$$
We denote by $H^{1}(\Omega)'$ the dual of $H^{1}(\Omega)$, and $H^{-1}(\Omega)$ the dual of $H^{1}_0(\Omega)$. For any $f\in H^{1}(\Omega)'$, we define a norm as
$$||f||_{H^1(\Omega)'}= \Big\{ \sup\,\langle f,u \rangle :u \in H^1(\Omega)\,,\,\,\,||u||_{H^1(\Omega)} \leq 1   \Big\} .$$
We introduce the solution space $W(0,T):= W_1(0,T)\times W_1(0,T)\times W_2(0,T)$ for the problem \eqref{maina}-\eqref{maine}, where 
\begin{align*}
 W_1(0,T)&= \Big\{w: \, w \in L^{2}(0,T;H^{1}(\Omega))\,, \,\, \dfrac{\partial w}{\partial t} \in L^{2}(0,T;(H^{1}(\Omega))^{\prime})  \Big\}\,, \\
 W_2(0,T)& = \Big\{w:  w \in L^{2}(0,T;H^{1}_0(\Omega))\,, \,\, \dfrac{\partial w}{\partial t} \in L^{2}(0,T;H^{-1}(\Omega))  \Big\}\,.
 \end{align*}
Note that, the space $W_i(0,T)~(i=1,2)$ is a Hilbert space for the graph norm, see~\cite{jllions1968}.
\begin{defi}[Weak solution]\label{defi:weak-solun}
\normalfont{
Let $I_0\in H^1(\Omega)$. We say that a triple $(I,u,v)$ is a weak solution of \eqref{maina}-\eqref{maine} if
\begin{itemize}
\item[(i)] $ I, u \in W_1(0,T)$, and $v\in W_2(0,T)$.
\item[(ii)] For all $\Psi \in H^1(\Omega)$ and $\Phi \in H_0^1(\Omega)$, there holds 
\begin{align}
\label{mainaweak}
&\int_0^T\big \langle \frac{\partial I}{\partial t},\Psi \big \rangle\, dt + \int_{\Omega_T} g(u)\,\nabla I\cdot\nabla \Psi \,dx\,dt + 2\int_{\Omega_T}  v \Psi \,dx\,dt=0\,,   \\
\label{mainbweak}
&\int_0^T \big \langle \frac{\partial u}{\partial t},\Psi \big \rangle\,dt  -\int_{\Omega_T} h(|\nabla I_{\xi}|^{2})\Psi \,dx\,dt  +\int_{\Omega_T} u\Psi \,dx\,dt \nonumber \\
& \hspace{5.5cm}+ \int_{\Omega_T} \nabla u \cdot \nabla \Psi \,dx\,dt  = 0\,, \\
\label{maincwaek}
&\int_0^T \big \langle \frac{\partial v}{\partial t}, \Phi \big \rangle\,dt + \int_{\Omega_T} \nabla v \cdot \nabla \Psi\, dx\,dt + \int_{\Omega_T} (I_0-I) \Psi\, dx\,dt=0\,.
\end{align}
\item[(iii)]  \eqref{maine} holds.
\end{itemize}
}
\end{defi}
\subsection{Discretized system and existence of its weak solution}\label{subsec:discretization}
To prove existence of a weak solution of \eqref{maina}-\eqref{maine}, we use semi discretization in time with a semi-implicit Euler method. Let $0=t_0<t_1<t_2 \dots <t_N =T$ be a uniform partition of $[0,T]$ with time-step size $\tau=\frac{T}{N}$ for some $N \in \mathbb{N}^+$, i.e, $t_n=\tau\,n$ for $0 \leq n \leq N$. For simplicity, in the following we write $L^2,$ $H^1,$ $H^1_0$ instead of $L^2(\Omega),$ $H^1(\Omega),$ $H^1_0(\Omega)$  respectively. We may then consider time discretization of \eqref{maina}-\eqref{maine} as follows: for $0\le n\le N-1\,,$ $\forall\, \Psi\in H^1$, and $ \forall\, \Phi \in  H^1_0$ iterate:
\begin{itemize}
\item [(i)] compute $u_{n+1} \in H^1$ such that
\begin{align}\label{eq:disc eq U}
\big(\frac{1}{\tau}(u_{n+1}-u_n),\Psi \big)_{L^2} + \big(\nabla u_{n+1},\nabla \Psi \big)_{L^2}+\big( u_{n+1},\Psi \big)_{L^2} \nonumber\\
 =\big( h(|\nabla G_{\xi}*I_{n}|^2),\Psi \big)_{L^2}\,,
\end{align}
\item[(ii)]compute $v_{n+1} \in H^1_0$ such that
\begin{align}\label{eq:disc eq V}
\big(\frac{1}{\tau}v_{n+1},\Phi \big)_{L^2} + \big(\nabla v_{n+1},\nabla \Phi \big)_{L^2}
 =\big(\frac{1}{\tau}v_{n} - (I_0-I_n),\Phi \big)_{L^2}\,
\end{align}
\item[(iii)] compute $I_{n+1} \in H^1$ 
\begin{align}\label{eq:disc eq I}
\big(\frac{1}{\tau}I_{n+1},\Psi \big)_{L^2} + \big( g(u_{n}) \nabla I_{n+1},\nabla \Psi \big)_{L^2} 
 =\big(\frac{1}{\tau}I_{n} - 2v_{n},\Psi \big)_{L^2}\,.
\end{align}
\end{itemize}
Thanks to \eqref{bound:h}, the existence of $u_{n+1}$ resp. $v_{n+1}$ in step ${\rm (i)}$ resp. in step ${\rm (ii)}$ easily follows from Lax-Milgram lemma. Moreover, these solutions are unique. To show the well-posedness
of \eqref{eq:disc eq I}, we use standard Lax-Milgram lemma. Define a continuous bilinear form $\mathcal{C}:H^{1} \times H^1 \rightarrow \mathbb{R} $, and a linear bounded functional $\ell:H^1\rightarrow\mathbb{R}$ as
\begin{align}
 \mathcal{C}(\Theta,\Psi) = \big(\frac{1}{\tau} \Theta ,\Psi \big)_{L^2} + \big(g(u_{n}) \nabla \Theta,\nabla \Psi \big)_{L^2}\,; \quad 
 {\ell}(\Psi)=\big(\frac{1}{\tau}I_{n}-2v_{n},\Psi \big)_{L^2}\,.
\end{align}
Then, by Lax-Milgram lemma, there exists a unique $I_{n+1}\in H^1$ satisfying \eqref{eq:disc eq I} provided $\mathcal{C}$ is $H^1$-coercive. Indeed, 
thanks to Gagliardo-Nirenberg inequality,there exists a constant $C>0$, depending only on $\Omega$ and $G_\xi$, such that
\begin{align}
 \dfrac{1}{  1+ \frac{C||u_{n}||_{H^1}}{k^2}  } \leq  g(u_{n})\,, \label{bound:g-lower}
\end{align}
and hence $\mathcal{C}(\Psi,\Psi)
\geq \frac{1}{\tau} ||\Psi||^2_{L^2}+   \dfrac{1}{  1+ \frac{C||u_{n}||_{H^1}}{k^2}  }||\nabla \Psi||^2_{L^2}
\geq  \widetilde{C}||\Psi||^2_{H^1} $
where $\widetilde{C}:= \min \Big\{ \frac{1}{\tau},  \dfrac{1}{  1+ \frac{C||u_{n}||_{H^1}}{k^2}  }\Big\}$ yielding the $H^1$-coercivity of $\mathcal{C}$. 
\subsection{A-priori estimates}
 Let $0 \leq n \leq N-1$.  Choose $\Psi = u_{n+1}$ as test function in \eqref{eq:disc eq U} and use the algebraic
identity 
\begin{align}\label{algebric Identity}
\langle a,a-b \rangle = \frac{1}{2}\big( |a|^2 - |b|^2 +|a-b|^2 \big)\,, \quad  \forall \,a,b \in \mathbb{R}^p\,(p\ge 1)\,,
\end{align}
Cauchy-Schwarz and Young's inequalities, \eqref{bound:h}, and then sum over $n=0,1,2, \dots ,\newline
j-1$ with $1 \leq j \leq N$. The result is
\begin{align}
||u_j||^2_{L^2}+\sum_{n = 0}^{j-1} ||u_{n+1}-u_n||^2_{L^2} + {2 \tau} \sum_{n = 0}^{j-1} ||\nabla u_{n+1}||^2_{L^2} \nonumber\\
 \le  ||u_0||^2_{L^2}  + |\Omega| \sum_{n=0}^{j-1} \tau \le ||u_0||^2_{L^2}  + T\,|\Omega|\,.
\end{align}
Hence, there exists a constant $C_1>0$, independent of $\tau >0$, such that 
\begin{align}\label{eq:estimate U final}
\max_{1\le j\le N} ||u_j||^2_{L^2}+\displaystyle\sum_{n = 0}^{N-1} ||u_{n+1}-u_n||^2_{L^2}+{2\tau} \sum_{n = 0}^{N-1} ||\nabla u_{n+1}||^2_{L^2} \leq C_1\,.
\end{align}
Similarly, by choosing~(formally) $\Psi =-\Delta u_{n+1}$ in \eqref{eq:disc eq U} and using the integration by parts formula, \eqref{algebric Identity}, 
Cauchy-Schwarz and Young's inequalities, \eqref{bound:h}, and then summing over $n=0,1,2, \ldots ,j-1$ with $1\le j\le N$,  we obtain

\begin{align}\label{eq: estimate for grad U}
\max_{1\le j\le N}||\nabla u_{j}||^2_{L^2}+ \sum_{n = 0}^{N-1} ||\nabla (u_{n+1}-u_n)||^2_{L^2} +{\tau} \sum_{n = 0}^{N-1} ||\Delta u_{n+1}||^2_{L^2} \leq C_2\,,
\end{align}
where $C_2>0$ is a constant, independent of $\tau>0$. 
In view of \eqref{eq:estimate U final} and \eqref{eq: estimate for grad U}, there exists a constant $C_3>0$, independent of $\tau$, such that
\begin{align}
\max_{0\le j\le N}\|u_j\|_{H^1(\Omega)}\le C_3\,. \label{bound:h1-norm-un}
\end{align}
Therefore, from \eqref{bound:g-lower}, the positive lower bound for $g(u_n)$ is given by
\begin{align}
 \nu:= \frac{1}{ 1 + \frac{CC_3}{k^2}} \le g(u_n)\quad (0\le n\le N)\,. \label{bound:lower-g-u_n}
\end{align}

Again, one can use the test function $\Phi = v_{n+1}$ and $\Phi =-\Delta v_{n+1}$~(formally) in \eqref{eq:disc eq V}, and proceed as above~(under the cosmetic changes) to obtain $(1 \leq j \leq N)$:
\begin{align}
&||v_j||^2_{L^2}+ \sum_{n = 0}^{j-1} ||v_{n+1}-v_n||^2_{L^2} + {\tau} \sum_{n = 0}^{j-1} ||\nabla v_{n+1}||^2_{L^2} \nonumber\\
&\hspace{4cm}\le \|v_0\|_{L^2}^2 + C_P^2 T \|I_0\|_{L^2}^2 + \tau C_P^2 \sum_{n=0}^{j-1}\|I_n\|_{L^2}^2\,, \label{eq: estimate V step 1} \\
&||\nabla v_{j}||^2_{L^2}+ \sum_{n = 0}^{j-1} ||\nabla (v_{n+1}-v_n)||^2_{L^2} +{\tau}\sum_{n = 0}^{j-1} ||\Delta v_{n+1}||^2_{L^2} \nonumber\\
&\hspace{5cm}\leq||\nabla v_{0}||^2_{L^2} + T \|I_0\|_{L^2}^2 +\tau \sum_{n = 0}^{j-1} ||I_n||^2_{L^2}\,, \label{eq: estimate gradV step 1}
\end{align}
where $C_P$ is the Poincar\'{e} constant.
Note that, because of coupled system, we are not able to find the bound of $v_n$'s from \eqref{eq: estimate V step 1} and \eqref{eq: estimate gradV step 1}. May be we need to combine with the estimates coming 
from \eqref{eq:disc eq I} by choosing appropriate test functions. Taking $\Psi = I_{n+1}$ in \eqref{eq:disc eq I} and using \eqref{algebric Identity}, \eqref{bound:lower-g-u_n}, and \eqref{eq: estimate V step 1},
we have 
\begin{align}
&||I_j||^2_{L^2}+ \sum_{n = 0}^{j-1} ||I_{n+1}-I_n||^2_{L^2} + {2\tau \nu} \sum_{n = 0}^{j-1} ||\nabla I_{n+1}||^2_{L^2}  \notag \\
& \leq ||I_0||^2_{L^2} + {2\tau}\sum_{n = 0}^{j-1} ||I_{n+1}||^2_{L^2} \nonumber\\
&\hspace{3cm}+ 2\tau\, \sum_{n=0}^{j-1}\Big(\|v_0\|_{L^2}^2 + C_P^2 T \|I_0\|_{L^2}^2 + \tau C_P^2 \sum_{k=0}^{n}\|I_{k+1}\|_{L^2}^2\Big) \notag \\
& \leq C_4 + C_5 \tau \sum_{k=0}^{j-1}\|I_{k+1}\|_{L^2}^2\,,
\end{align}
for some constants $C_4, C_5>0$, only depend on $T, C_P, v_0, I_0$. We now apply discrete Gronwall's lemma~(implicit form) and obtain the following: there exists $\tau_1>0$ such that for all time step
sizes $0<\tau<\tau_1$, 
\begin{align}
\max_{0\le j\le N}||I_j||^2_{L^2}+ \sum_{n = 0}^{N-1} ||I_{n+1}-I_n||^2_{L^2} + {2\tau \nu} \sum_{n = 0}^{N-1} ||\nabla I_{n+1}||^2_{L^2} \leq  C_6\,, \label{eq: estimate I final}
\end{align}
where $ C_6 > 0 $ is independent of $\tau > 0$.  We now use \eqref{eq: estimate I final} in \eqref{eq: estimate V step 1} and \eqref{eq: estimate gradV step 1},
and get
\begin{align}
\max_{0\le j\le N} ||v_j||^2_{H^1}+\sum_{n = 0}^{N-1} ||v_{n+1}-v_n||^2_{H^1}+ \tau \sum_{n = 0}^{N-1} \Big(||\nabla v_{n+1}||^2_{L^2} + \|\Delta v_{n+1}\|_{L^2}^2\Big) \leq C_7\,, \label{eq:estimate V final}
\end{align}
for some constant $C_7>0$, independent of $\tau$ with $0<\tau < \tau_1$.
Again one may choose $\Psi =-\Delta I_{n+1}$~(formally) in \eqref{eq:disc eq I} and use the integration by parts formula, \eqref{algebric Identity}, Cauchy-Schwarz and Young's inequalities,  \eqref{bound:lower-g-u_n},
\eqref{bound:h1-norm-un}, \eqref{eq:estimate V final}  and the fact that $\|\nabla g(u_n)\|_{L^\infty} \le C\big(\xi, \Omega, k\big)\|u_n\|_{H^1}$ to get, after the application of discrete Gronwall's 
lemma (implicit form), 

\begin{align}
\max_{1\le j\le N}||\nabla I_{j}||^2_{L^2}+ \sum_{n = 0}^{N-1} ||\nabla (I_{n+1}-I_n)||^2_{L^2} + \tau \sum_{n = 0}^{j-1}||\Delta I_{n+1}||^2_{L^2}  \leq C_8\,, \label{eq: estimate for grad I}
\end{align}
with $ C_8 > 0 $, independent of $ \tau $ with $0<\tau<\tau_2$ for some $\tau_2>0$. Note that, for $\tau_0=\min\{ \tau_1, \tau_2\}$, \eqref{eq: estimate I final}, \eqref{eq:estimate V final} and 
\eqref{eq: estimate for grad I} all hold true. Putting things together, we arrive at the following theorem.
\begin{thm}\label{thm:a-priori}
\normalfont{
 Let $I_0\in H^1$. Then there exist constants $\tau_0>0$ and $C>0$ such that for all time step sizes $ \tau>0 $ with $0<\tau<\tau_0$, there holds
 \begin{align}
 &  \max_{1\le j\le N} \Big( \|u_j\|_{H^1}^2 + \|v_j\|_{H^1}^2 + \|I_j\|_{H^1}^2 \Big) \nonumber\\
 & \hspace{2cm} + \sum_{n=0}^{N-1} \Big( ||u_{n+1}-u_n||^2_{H^1} + ||v_{n+1}-v_n||^2_{H^1} + ||I_{n+1}-I_n||^2_{H^1}\Big) \notag \\
  &\hspace{2cm} + \tau \sum_{n=0}^{N-1} \Big(\|\Delta u_{n+1}\|_{L^2}^2 + \|\Delta v_{n+1}\|_{L^2}^2 + \|\Delta I_{n+1}\|_{L^2}^2\Big) \le C\,, 
 \end{align}
where $u_j$, $v_j$ and $I_j$ solves \eqref{eq:disc eq U}, \eqref{eq:disc eq V} and \eqref{eq:disc eq I} respectively. 
}
\end{thm}
\subsection{Continufication and existence of weak solution} Let $\{t_n\}_{n=0}^N$ be a uniform partition of $[0,T]$ with time-step size $\tau$ as described in Subsection \ref{subsec:discretization}. 
 For any sequence $\{\mathrm{x}_n\} \subset \mathbb{X}$, where $\mathbb{X}$ is a Banach space, we define the difference quotient $d_t\mathrm{x}_{n+1}=\dfrac{\mathrm{x}_{n+1}-\mathrm{x}_n}{\tau}$
 for $0 \leq n \leq N-1$.  The globally time interpolant $X_\tau \in C(\mathbb{X})$ of $\{\mathrm{x}_n\}$ is defined via
\begin{align*}
\mathrm{X}_\tau(t):=\dfrac{t-t_n}{\tau} {\mathrm{x}}_{n+1}+\dfrac{t_{n+1}-t}{\tau} {\mathrm{x}}_{n}\,, \quad \forall\, t \in (t_n,t_{n+1}]\,.
\end{align*}
Moreover, we define the piecewise constant in time interpolants ${\mathrm{X}}^{+}_{\tau}(t)$ and ${\mathrm{X}}^{-}_{\tau}(t)$ as follows:
\begin{align*}
{\mathrm{X}}^{+}_{\tau}(t):={\mathrm{x}}_{n+1}\,, \quad  {\mathrm{X}}^{-}_{\tau}(t):={\mathrm{x}}_n\,, \quad \forall\, t \in (t_n,t_{n+1}]\,,
\end{align*}
with ${\mathrm{x}}_{-1}={\mathrm{x}}_0$ and ${\mathrm{x}}_{N+1}=0$.
\vspace{.1cm}
With the above notation, we now show the boundedness of the sequence $\{ \partial_t \mathcal{U}_\tau \}_{\tau>0}$ in $L^2(H^{1}(\Omega)')$. To do so,  let $\Psi \in H^{1}(\Omega)'\setminus \{0\}$, and 
$ t \in (t_n,t_{n+1}]$ where $0\le n\le N-1$. Then from \eqref{eq:disc eq U}, we have, by using Cauchy-Schwarz inequality and \eqref{bound:h},
\begin{align}
 \big \langle \partial_t \mathcal{U}_{\tau}, \Psi \big \rangle 
 & \le \Big( |\Omega|^{\frac{1}{2}} + \max_{0\le n\le N-1}\|u_{n+1}\|_{H^1}\Big) \|\Psi\|_{H^1}\,,
\end{align}
and hence, thanks to  Theorem \ref{thm:a-priori}, 
\begin{align}
|| \partial_t \mathcal{U}_{\tau}||^2_{L^2({H^{1}(\Omega)'})} \leq & C\Big( \tau \sum_{n = 0}^{N-1} |\Omega| + T  \max_{0\le n\le N-1}\|u_{n+1}\|^2_{H^1}\Big)\le C\,. \label{eq: estimate of disc dU}
\end{align}
Similarly, the sequences  $\{ \partial_t \mathcal{V}_\tau \}_{\tau>0}$ and  $\{ \partial_t \mathcal{I}_\tau \}_{\tau>0}$ are bounded 
in $L^2(H^{-1}(\Omega))$ and $L^2(H^{1}(\Omega)')$ respectively. Moreover, there exists a constant $C>0$ such that for all $0< \tau < \tau_0$,
\begin{align}
 || \partial_t \mathcal{V}_{\tau}||^2_{L^2({H^{-1}(\Omega)})} + || \partial_t \mathcal{I}_{\tau}||^2_{L^2({H^{1}(\Omega)'})} \le C\,.\label{eq: estimate of disc dV+dI}
\end{align}

Thanks to \eqref{eq: estimate of disc dU}, \eqref{eq: estimate of disc dV+dI} and Theorem~\ref{thm:a-priori}, by using classical results of compact inclusion in Sobolev spaces \cite{raadams1975},
there exists $(I,u,v)\in W_1(0,T)\times W_1(0,T)\times W_2(0,T)$ such that along a
subsequence (still we denote it by the same index) the following hold:
\begin{align}\label{eq:convergence-continufication}
\begin{cases}
\mathcal{I}_{\tau}^{+},\mathcal{I}_{\tau}^{-}\rightharpoonup^{ \ast} I \quad  \text{in}\,\,\, L^{\infty}(H^1)\,; 
\quad  \partial_t \mathcal{I}_{\tau} \rightharpoonup \partial_t I \quad  \text{in}\,\,\, L^2({H^{1}(\Omega)'})\,, \\
\mathcal{I}_{\tau},\mathcal{I}_{\tau}^{\pm}\rightharpoonup I \quad  \text{in}\,\,\, L^{2}(H^1)\,; \quad 
\mathcal{I}_{\tau},\mathcal{I}_{\tau}^{\pm} \rightarrow I \quad  \text{in}\,\,\, L^{2}(L^2)\,, \\
|\nabla G_{\xi} \ast \mathcal{I}^{-}_{\tau}|^2 \rightarrow |\nabla G_{\xi} \ast I|^2 \quad  \text{in}\,\,\,L^{2}(L^2) \,\,\text{and a.e.}\, (0,T)\times \Omega\,,\\
h(|\nabla G_{\xi} \ast \mathcal{I}_{\tau}|^2) \rightarrow h(|\nabla G_{\xi} \ast I|^2)\quad  \text{in}\,\,\, L^{2}(L^2)\,, \\
\mathcal{U}_{\tau}^{+},\mathcal{U}_{\tau}^{-}\rightharpoonup^{ \ast} u\quad  \text{in}\,\,\, L^{\infty}(H^1)\,; \quad 
\partial_t \mathcal{U}_{\tau} \rightharpoonup \partial_t u \quad  \text{in}\,\,\, L^2({H^{1}(\Omega)'})\,, \qquad (\tau \rightarrow 0)\,, \\
\mathcal{U}_{\tau},\mathcal{U}_{\tau}^{\pm}\rightharpoonup u \quad  \text{in}\,\,\, L^{2}(H^1)\,; \quad 
\mathcal{U}_{\tau},\mathcal{U}_{\tau}^{\pm} \rightarrow u \quad  \text{in}\,\,\, L^{2}(L^2)\,, \\
g(\mathcal{U}^{+}_{\tau}) \rightarrow g(u)  \quad  \text{in}\,\,\, L^{2}(L^2)\,, \\
\mathcal{V}_{\tau}^{+},\mathcal{V}_{\tau}^{-}\rightharpoonup^{ \ast} v \quad  \text{in}\,\,\, L^{\infty}(H^1_0)\,;
\quad \partial_t \mathcal{V}_{\tau} \rightharpoonup \partial_t v \quad  \text{in}\,\,\, L^2({H^{-1}(\Omega)})\,, \\
\mathcal{V}_{\tau},\mathcal{V}_{\tau}^{\pm}\rightharpoonup v \quad  \text{in}\,\,\, L^{2}(H^1_0)\,; 
\quad \mathcal{V}_{\tau},\mathcal{V}_{\tau}^{\pm} \rightarrow v \quad  \text{in}\,\,\, L^{2}(L^2)\,.
\end{cases}
\end{align}
For any $\Psi \in H^1$ and $\Phi \in H_0^1$, we rewrite \eqref{eq:disc eq U}, \eqref{eq:disc eq V} and \eqref{eq:disc eq I} in terms of $\mathcal{U}_\tau$, ${\mathcal{U}}_{\tau}^{+}$, 
${\mathcal{U}}_{\tau}^{-}$,  $\mathcal{V}_\tau$, ${\mathcal{V}}_{\tau}^{+}$, ${\mathcal{V}}_{\tau}^{-}$,  $\mathcal{I}_\tau$, ${\mathcal{I}}_{\tau}^{+}$, and ${\mathcal{I}}_{\tau}^{-}$, integrate with respect to
the time variable and make use of \eqref{eq:convergence-continufication} to pass to the limit as $\tau \rightarrow 0$  in the resulting variational formulation to arrive at 
\eqref{mainaweak}, \eqref{mainbweak} and \eqref{maincwaek}. In other words, the triple $(I,u,v)$ is a weak solution of \eqref{maina}-\eqref{maine}.
\subsection{Uniqueness of weak solution} 
We use standard methodology \cite{lcevans1998} to prove the uniqueness of weak solution of \eqref{maina}-\eqref{maine}. Let $(I_1,u_1,v_1)$ and $(I_2,u_2,v_2)$ be two sets of solution for the system \eqref{maina}-\eqref{maine}
with $I_1\neq I_2$, $u_1\neq u_2$ and $v_1\neq v_2$. Let $
\tilde{I}:=I_1-I_2\,,\, \widetilde{U}:=u_1-u_2\,,\, \widetilde{V}:=v_1-v_2$ and $ g_i=g(u_i)~(i=1,2)$. 
Then, the following equations hold in the sense of distribution:
\begin{align}
&\dfrac{\partial \tilde{I} }{\partial t} - \nabla(g_1\nabla \tilde{I}) =  \nabla((g_1-g_2)\nabla I_2)-2\widetilde{V}\,,  \hspace{2.7cm} \text{in}\,\,\, \Omega_T\,, \label{uni_maina_tilde} \\
&\dfrac{\partial \widetilde{U}}{\partial t} -\Delta\widetilde{U} = \Big(h(|\nabla G_{\xi}\ast I_1|^{2})-h(|\nabla G_{\xi}\ast I_2|^{2}) \Big)
-\widetilde{U}\,, \hspace{0.9cm} \text{in} \,\,\, \Omega_T \, , \label{uni_mainb_tilde}\\
&\dfrac{\partial\widetilde{V}}{\partial t}-\Delta \widetilde{V} =\tilde{I}\,,  \hspace{6.7cm} \text{in} \,\,\, \Omega_T\,, \label{uni_mainc_tilde} \\
& \tilde{I}(0,x)=0=\widetilde{U}(0,x)=\widetilde{V}(0,x)\,, \hspace{4.3cm} \text{in}\,\,\, \Omega\,. \label{ini_cond_tilde}
\end{align}
Note that
\begin{align*}
\nu_1:=\frac{1}{1 + \frac{C\|u_1\|_{L^\infty(H^1)}}{k^2}}\le g(u_1)\,; \quad || (g_1-g_2)||_{L^\infty} \leq C_{\xi}|| (u_1-u_2)||_{L^2}\,.
\end{align*}
Similar to a-priori estimates, we have for a.e. $t\in (0,T)$ with $\tilde{C}_I:= \|\nabla I\|_{L^\infty(L^2)}^2$,
\begin{align}
 \frac{1}{2}\frac{d}{dt}\|\tilde{I}\|_{L^2}^2 
 & \le \tilde{C}_I\,C(\nu_1) || (g_1-g_2)||_{L^\infty}^2  + \| \tilde{I}\|_{L^2}^2 + \|\widetilde{V}\|_{L^2}^2 \nonumber\\
 & \le  \widetilde{C}||\widetilde{U}||^2_{L^2}+||\tilde{I}||^2_{L^2}+||\widetilde{V}||^2_{L^2}\,.\label{eq:relation for I}
\end{align}
Similarly, by multiplying \eqref{uni_mainb_tilde} by $\widetilde{U}$ and integrating over $\Omega$, we get
\begin{align}
\frac{d}{dt}||\widetilde{U}||^2_{L^2} +2 ||\nabla \widetilde{U}||^2_{L^2} \le C ||h(|\nabla G_{\xi}\ast I_1|^{2})-h(|\nabla G_{\xi}\ast I_2|^{2})||_{L^2}^2 +  ||\widetilde{U}||^2_{L^2}\,.
\end{align}
By using Lipschitz continuity of $h$ along with Gagliardo-Nirenberg inequality, we have 
\begin{align}
 ||h(|\nabla G_{\xi}\ast I_1|^{2})-h(|\nabla G_{\xi}\ast I_2|^{2})||_{L^2}^2
 & \le C(c_h, \xi, I_1, I_2) \|\tilde{I}\|_{L^2}^2\,,
\end{align}
and hence 
\begin{align}
 \frac{d}{dt}||\widetilde{U}||^2_{L^2} \le C ||\tilde{I}||^2_{L^2}\,.\label{eq:relation for U}
\end{align}
In a similar way, one can easily deduce that 
\begin{align}
 \frac{d}{dt} || \widetilde{V}||^2_{L^2} \leq C\Big( || \tilde{I}||^2_{L^2} +  || \widetilde{V}||^2_{L^2}\Big)\,. \label{eq:relation for V}
\end{align}
We now add \eqref{eq:relation for I}-\eqref{eq:relation for V} and then apply Gronwall's inequality to infer that  $\tilde{I}\equiv 0$, $\widetilde{U}\equiv 0$ and $\widetilde{V}\equiv 0$. In other words,
weak solution of the proposed model \eqref{maina}-\eqref{maine} is unique. 

\section{Numerical Approximation}
\label{sec:numerical}
In this section, we discuss the numerical implementation of the proposed model \eqref{maina}-\eqref{maine}. To solve the model numerically, we have applied the generalized weighted average finite difference scheme which has a combined nature of the forward Euler method and the backward Euler method, at $n^{th}$ and $(n+1)^{th}$ step, respectively. The image variable $I$, edge variable $u$ and fidelity term $v$ are calculated at each iteration.

In the derivation of finite difference formulas, $h$ and $\tau$ are considered as spatial step size and time step size, respectively. Let $I(x_i,y_j,t_n)$ express the gray level of the image plane $I_{i,j}^n$, where $x_i=ih,y_j=jh$ and $t_n=n\tau$. Also, we consider ${v}_{i,j}^0=0$ and ${u}_{i,j}^0=G_{\xi}\ast|\nabla I^{0}_{i,j}|^{2}$. Derivatives are approximated by central difference formula as,
\begin{align*}
\displaystyle\frac{\partial I}{\partial x}& \approx \displaystyle\frac{I_{i+1,j}^n-I_{i-1,j}^n}{2h}\,,
\displaystyle\frac{\partial I}{\partial y} \approx \displaystyle\frac{I_{i,j+1}^n-I_{i,j-1}^n}{2h},\\
\displaystyle\frac{\partial^2 I}{\partial^2 x}& \approx \displaystyle\frac{I_{i+1,j}^n-2I_{i,j}^n+I_{i-1,j}^n}{h^2}\,,
\displaystyle\frac{\partial^2 I}{\partial^2 y} \approx \displaystyle\frac{I_{i,j+1}^n-2I_{i,j}^n+I_{i,j-1}^n} {h^2}.
\end{align*}
The discrete form of \eqref{maina}-\eqref{maine} can be written as;
\begin{equation}\label{eq:C-N}
\displaystyle
I_{i,j}^{n+1}-\frac{\tau}{2} \left((\nabla(g(u)\nabla I)^{n+1}_{i,j}-2\lambda v^{n+1}_{i,j}\right) = I_{i,j}^{n}+\frac{\tau}{2}\left((\nabla(g(u)\nabla I)^{n}_{i,j}-2\lambda v^{n}_{i,j}\right),
\end{equation}
\begin{equation}
\displaystyle
\frac{u_{i,j}^{n+1}-u_{i,j}^n}{\tau}=\varphi \left(  h(|\nabla {I^{n}_{\xi_{i,j}}}|^{2})-u_{i,j}^{n+1}+\frac{\psi^2}{2}\Delta u_{i,j}^{n+1}\right), 
\end{equation}
\begin{equation}
\frac{v_{i,j}^{n+1}-v_{i,j}^n}{\tau}=\Delta v_{i,j}^{n+1}-\left( I_{i,j}^n-I_{i,j}^{0}\right),
\end{equation}
where,
\begin{align}
\Delta w_{i,j}^{n}=\displaystyle\frac{w_{i+1,j}^n-2w_{i,j}^n+w_{i-1,j}^n}{h^2}+\displaystyle\frac{w_{i,j+1}^n-2w_{i,j}^n+w_{i,j-1}^n} {h^2},
\end{align}
for $w=u,v.$	
Moreover the diffusion term in \eqref{maina}-\eqref{maine} approximated using the central difference scheme as,
\begin{multline}
\nabla(g(u)\nabla I)_{i,j} = \frac{1}{h^2}( g_{i+\frac{1}{2},j}(I_{i+1,j}-I_{i,j})-g_{i-\frac{1}{2},j}(I_{i,j}-I_{i-1,j})+\\
g_{i,j+\frac{1}{2}}(I_{i,j+1}-I_{i,j})-g_{i,j-\frac{1}{2}}(I_{i,j}-I_{i,j-1}))\,, 
\end{multline}
with the Neumann boundary conditions: 
$$I_{i,0}-I_{i,1}=0, \quad I_{i,N-1}-I_{i,N}=0, \quad \text{for} \quad 0\leq i\leq M,$$
$$I_{0,j}-I_{1,j}=0, \quad I_{M-1,j}-I_{M,j}=0, \quad \text{for} \quad 0\leq j\leq N.$$
Further, to solve the algebraic system of the the form $AI^{n+1}=B$, generated from numerical discretization, hybrid Bi-Conjugate Gradient Stabilized solver has been used \cite{jain2015Elsevier,jain2015comparative}.\newline
Apart from the numerical discretization of \eqref{maina}-\eqref{maine}, a stopping criterion is needed to terminate the diffusion process. To achieve this goal, we start with an initial image $I_0$ and apply the system \eqref{maina}-\eqref{maine} repeatedly. This results in a family of smoother images ${I(t, x)};t>0$, which depicts refined versions of $I_0$. After sufficient iterations, changes between two consecutive iterations become redundant. At this point, the convergence of the iterative procedure has been achieved. Now, we use the following measure as a stopping criterion,
\begin{equation}\label{eq:4stopping}
\frac{{||I^{k+1}-I^k||}_2^2}{{||I^k||}_2^2}\leq \varepsilon\,.
\end{equation} 
Here, we use $\varepsilon = 10^{-4}$ as a fixed threshold. $I^k$ and $I^{(k+1)}$ depict the image planes at the $k^{th}$ and ${(k+1)}^{th}$ iteration, respectively.
\begin{algorithm}
\small
\begin{algorithmic}[1]
\STATE Inputs: $I^{(0)}=I_0$ as initial solution or initial noisy image,\\
	           $v^{(0)}=0$ as initial edge variable and $u^{(0)}=G_{\xi}\ast|\nabla I^{(0)}|^{2}$,
	           $\tau$ as time step, $k=0$.
	           
\STATE Result: Denoised image (I)
 
\STATE Calculate $I^{(k+1)}$, $u^{(k+1)}$ and $v^{(k+1)}$ from
\[
\begin{bmatrix}
    1-\frac{\tau}{2}A(I^{k+1}) & 0 & \tau \lambda \\
    0 & 1+\tau(1-\frac{\psi^2}{2}\Delta) & 0 \\
    0 & 0 & 1-\tau \Delta
  \end{bmatrix}
\begin{bmatrix}
    I^{k+1} \\
    u^{k+1} \\
   v^{k+1}
  \end{bmatrix}
  =
  \begin{bmatrix}
    I^{k}+\frac{\tau}{2}(A(I^{k})-2\lambda v^{k}) \\
    u^{k}+\tau \varphi h(||\nabla {I^k_\xi}||^{2})\\
    v^{k}-\tau(I^{k}-I^{0})
  \end{bmatrix}
\].
\STATE Solve the algebraic system (obtained using step 3), using Hybrid BiCGStab solver .
\STATE Check if 
		$$\frac{{||I^{k+1}-I^k||}_2^2}{{||I^k||}_2^2}\leq \varepsilon = 10^{-4}$$
then stop or go to Step 3.
\STATE end
\end{algorithmic}\caption{\textit{Algorithm of CPDE model for Image Denoising}}
\label{Algo:Algorithm for Hybrid BiCGStab Method}
\end{algorithm}
Parameter values used for the numerical experiments are mentioned in the caption of each figure. For the TV  and NLM  model, we have reproduced the results by utilizing the fact as mention in \cite{rudin1992nonlinear} and \cite{buades2005non} respectively. Apart from the parameters displayed in the captions, we have chosen $\varphi =1$ and $\xi=1$ for the present model and a uniform time step size $\tau =0.1$ for the current as well as the other discussed approaches.

\section{Results and Discussion}
\label{sec:result}
\begin{figure}[]
        \centering 
        \begin{subfigure}[b]{0.3\textwidth}           
                \includegraphics[scale=0.21]{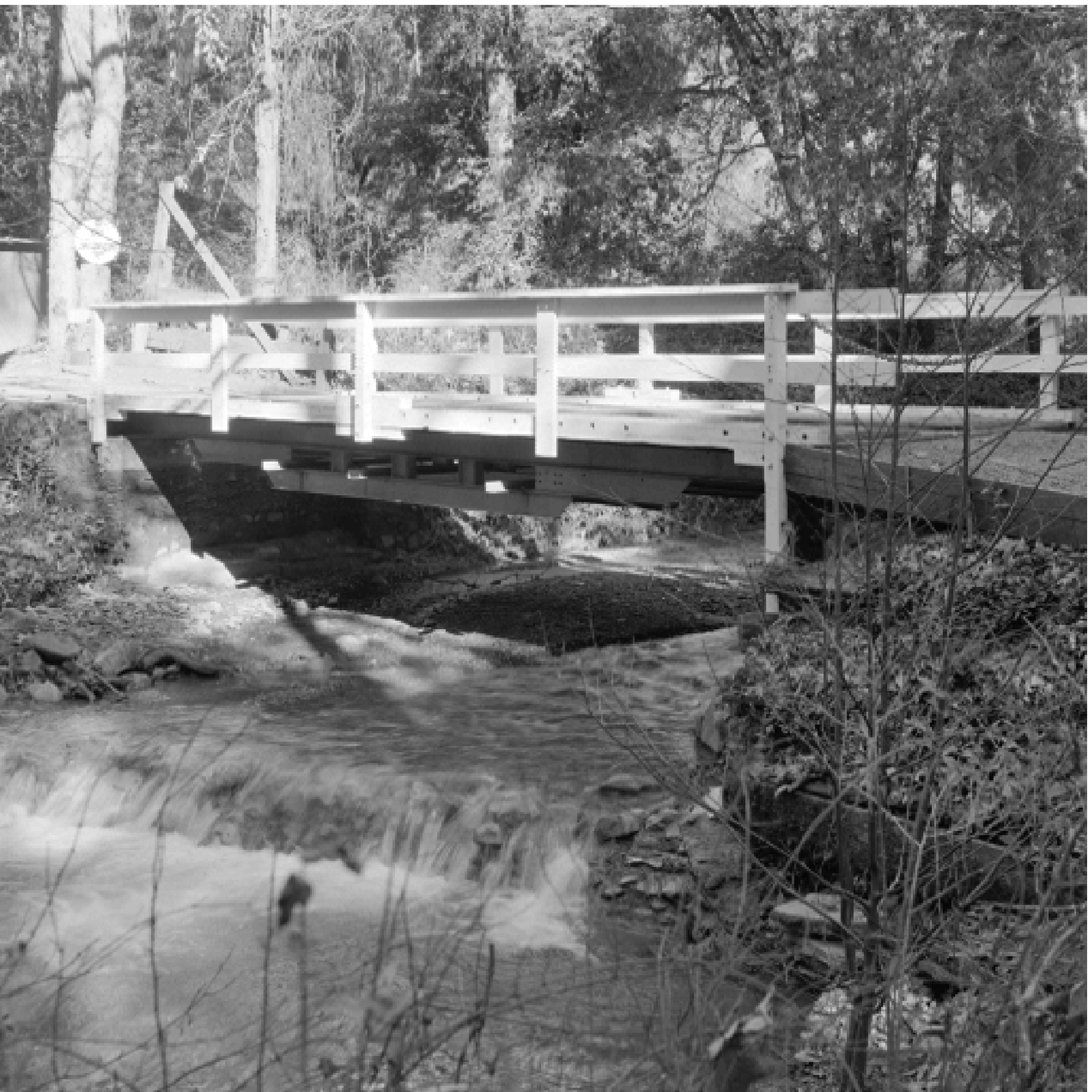}               
                \caption{Clear}
                \label{fig:3a}
        \end{subfigure}%
       \begin{subfigure}[b]{0.3\textwidth}
                \includegraphics[scale=0.21]{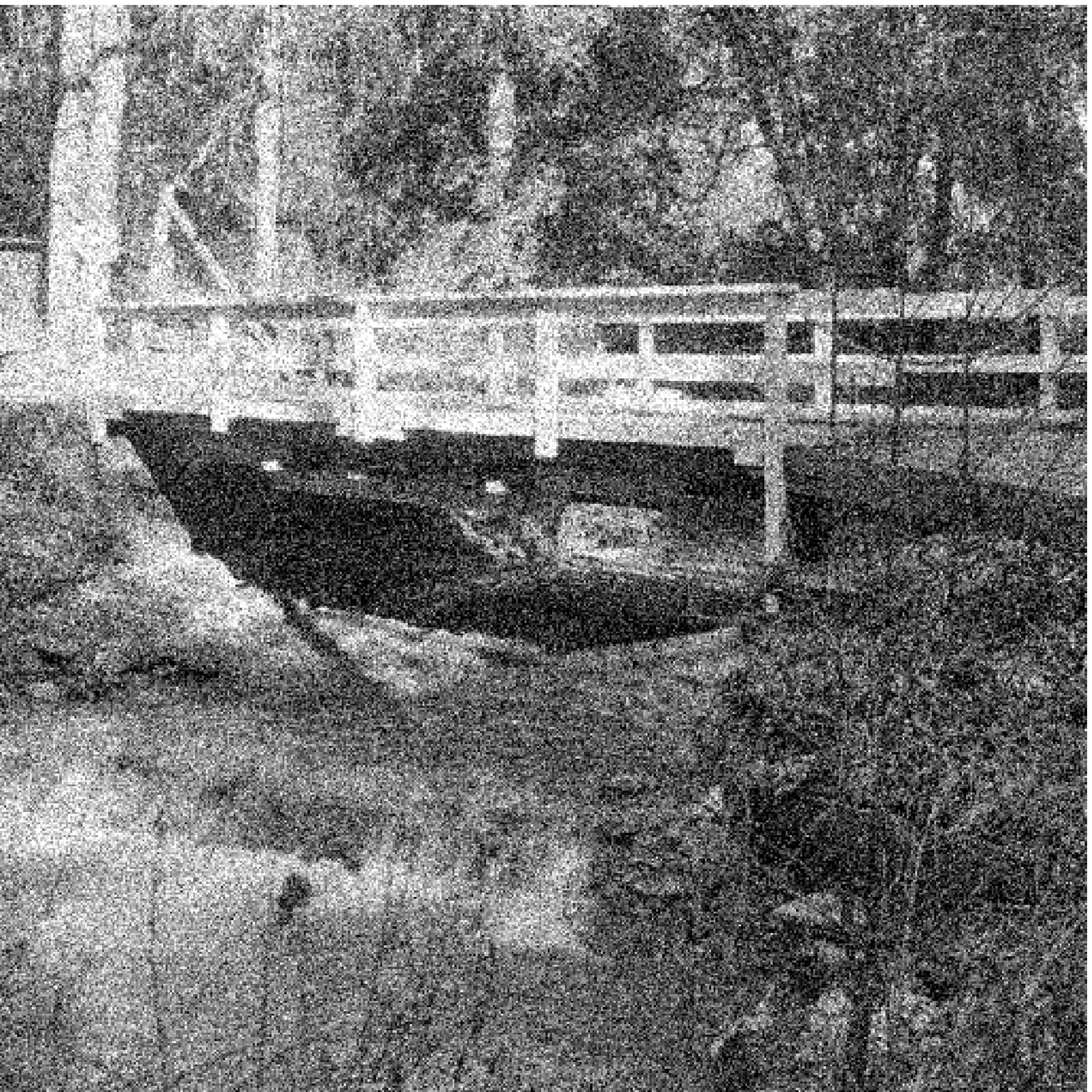}
                \caption{Noisy}
                \label{fig:3b}
        \end{subfigure}

        \begin{subfigure}[b]{0.3\textwidth}
                \includegraphics[scale=0.21]{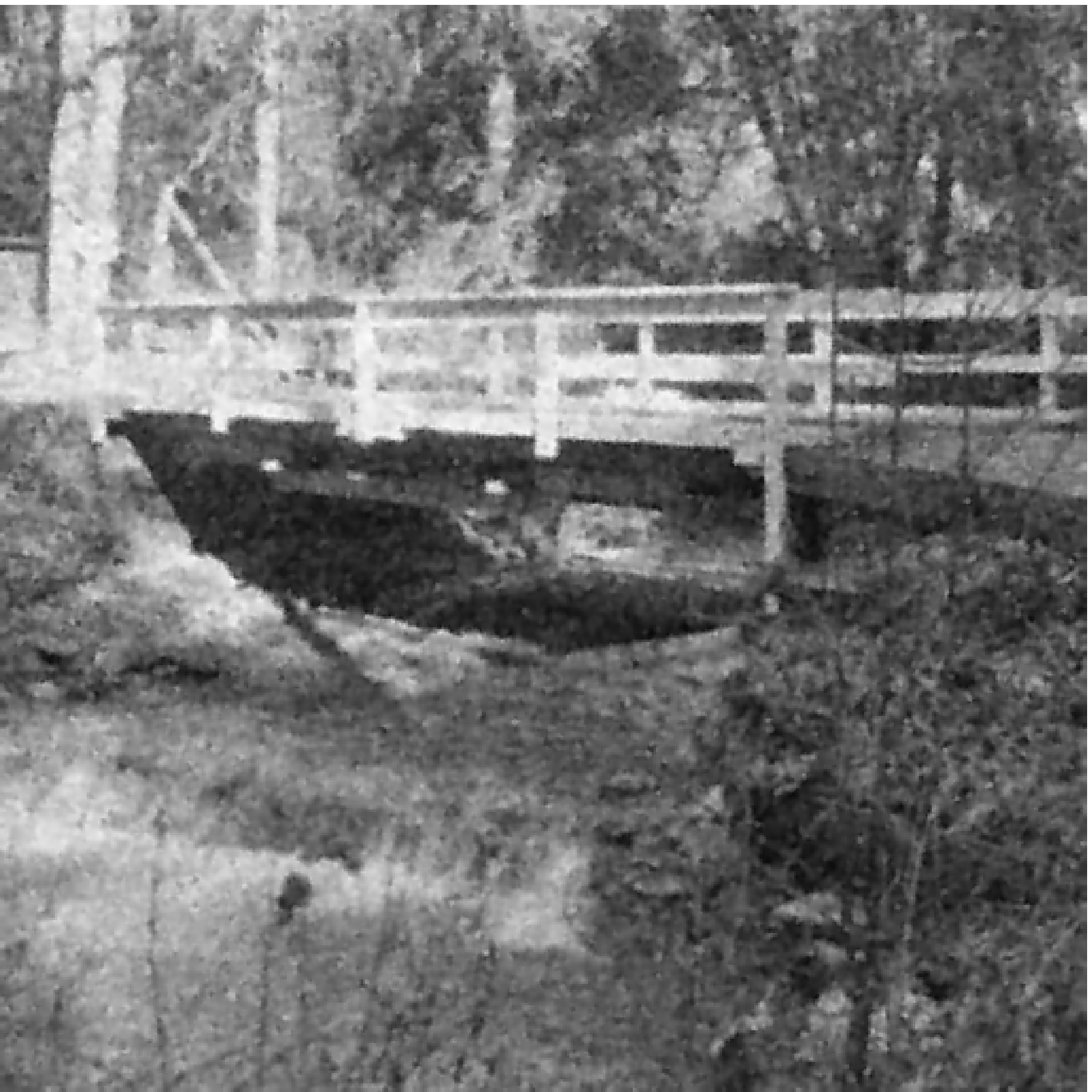}
                \caption{PV}
                \label{fig:3h}
        \end{subfigure}
        \begin{subfigure}[b]{0.3\textwidth}
                \includegraphics[scale=0.21]{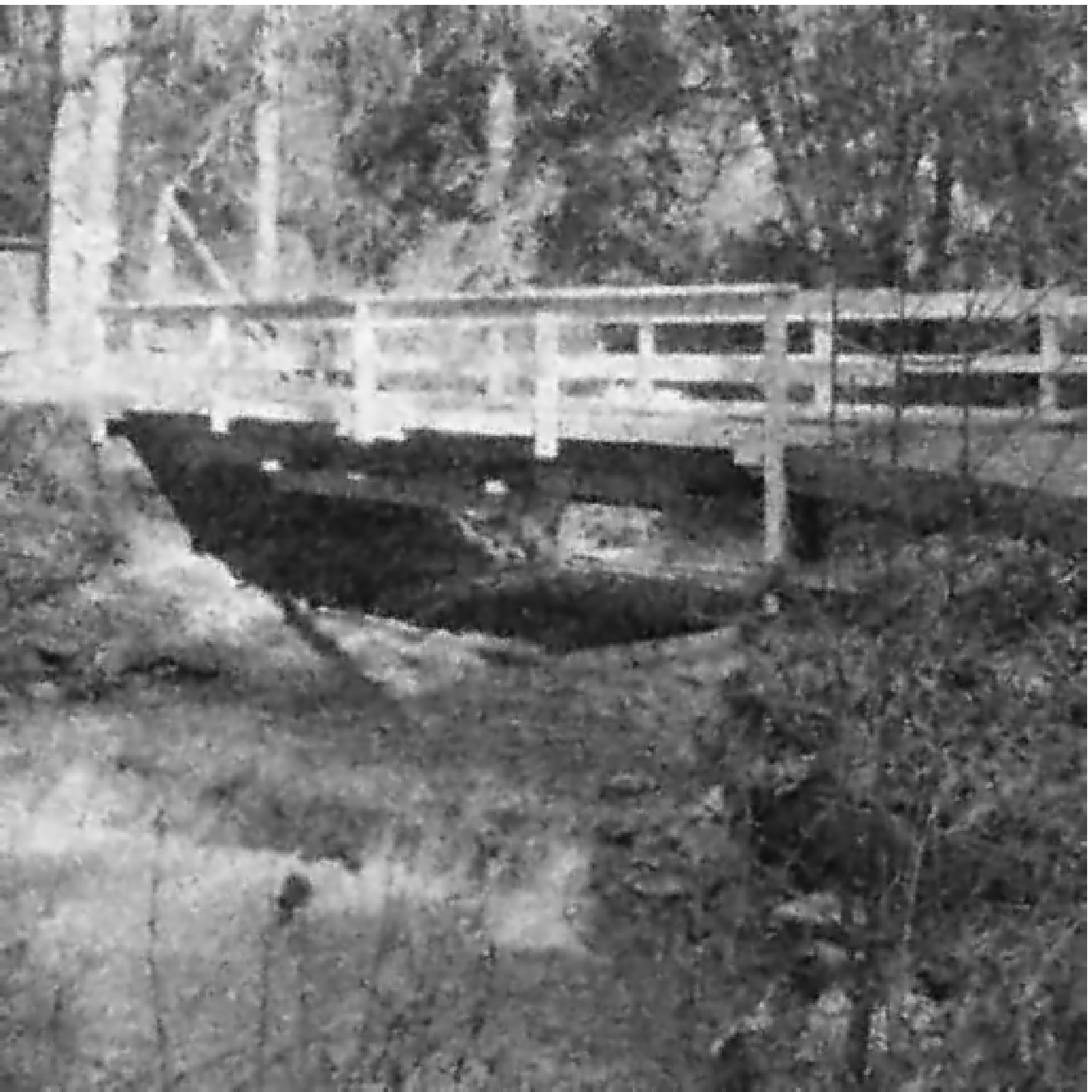}
                \caption{SYS}
                \label{fig:3i}
        \end{subfigure}
        \begin{subfigure}[b]{0.3\textwidth}
                \includegraphics[scale=0.21]{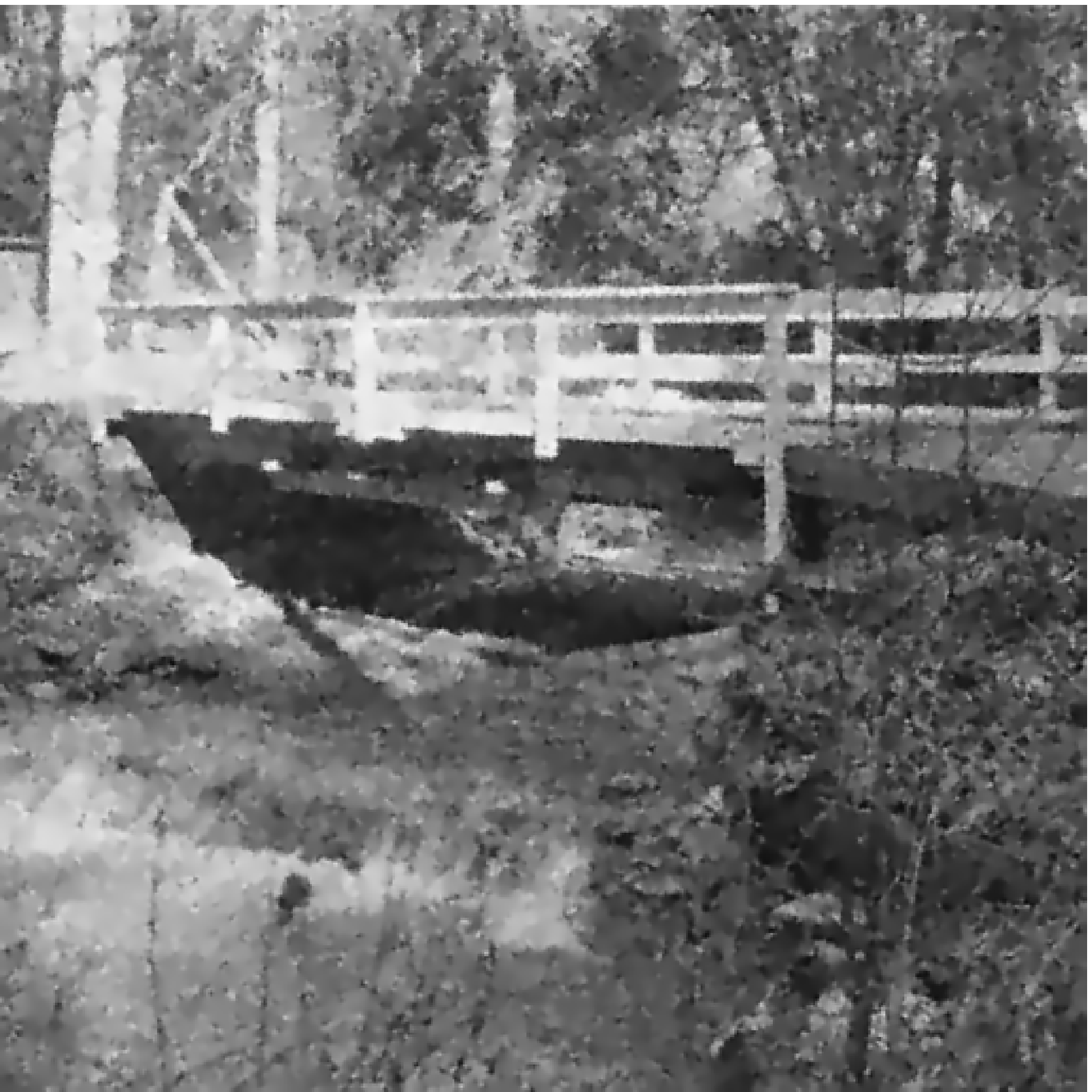}
                \caption{Proposed}
                \label{fig:3j}
        \end{subfigure}
        \caption{(a) Clear image; (b) Noisy image with Gaussian noise of mean $0.0$ and $\sigma=40$; Denoised image using (c) PV Model, $\lambda=0.01, K=8$; (d) SYS Model, $\lambda=4, K=8$; (e) Proposed model, $\psi=0.1, k=4.45.$}\label{fig:walkbridge_40}
\end{figure}

The proposed CPDE model was employed using the finite difference scheme given in the previous section. Image denoising using \eqref{maina}-\eqref{maine} was compared with the results of other state-of-art models available in the literature. Especially, TV model \cite{rudin1992nonlinear}, NS model \cite{nitzberg1992nonlinear}, Luo model \cite{luo2006coupled}, RD model \cite{guo2011reaction}, PV model \cite{prasath2014system}, SYS model \cite{sun2016class} are considered for the comparison. Apart from the diffusion-based approaches, non-iterative approach (Non-Local Means method) \cite{buades2005non} is also used for comparison with the present method. Since the space-time regularization based proposed model is claimed to be an improvement over the existing CPDE models, our main aim was to compare the edge detection and denoising results with these models (NS Model, and Luo model). In this process, the considered non-linear diffusion models are solved by the existing numerical schemes. And also, to stop the iterative process of each smoothing algorithm discussed stopping criteria equation \eqref{eq:4stopping} is employed. The effectiveness of results was evaluated through several standard gray level test images which are degraded with additive Gaussian noise of zero mean and different levels of standard deviations. We have artificially added additive Gaussian noise of different standard deviations ranging from $20$ to $50$ by using our MATLAB program.
Also, to evaluate the efficiency of our model, quantitative comparisons in terms of Peak signal to noise ratio (PSNR)\cite{gonzalez2002digital}, mean structural similarity index measure (MSSIM)\cite{wang2004image}, Gradient of peak signal to noise ratio ($\text{PSNR}_\text{Grad}$)\cite{tai2006image}, Improvement in signal to noise ratio (ISNR)\cite{gonzalez2002digital} have been shown with existing models. A higher value of quantity metrics suggests that the filtered image is closer to the noise-free image.\\
(a). Peak signal to noise ratio (PSNR)\cite{gonzalez2002digital} can measure the match between the clean ($I$) and denoised image ($\hat{I}$),
\begin{align*}
\text{PSNR}(I,\hat{I}) = 10\text{log}_{10}\frac{\text{MN}|\text{max}(I)-\text{min}(I)|^2}{{||\hat{I}-I||}^2_{L^2}}\,.
\end{align*}
(b). In addition to PSNR, gradient of peak signal to noise ratio ($\text{PSNR}_\text{Grad}$) is also used to measure the match between the derivatives of reconstructed and true image, defined as
\begin{align*}
\text{PSNR}_{\text{Grad}}=\frac{1}{2}(\text{PSNR}(I_x,\hat{I}_x)+\text{PSNR}(I_y,\hat{I}_y))\,,
\end{align*}
where $(I_x,I_y)$ and $(\hat{I}_x,\hat{I}_y)$ are the derivatives of ground truth image $(I)$ and denoised image $(\hat{I})$.\\
(c). Structural similarity index (SSIM) \cite{wang2004image}, is used to calculate the similarity between structure of clean and reconstructed images, and can be given as,
\begin{align*}
\text{SSIM}(X,Y)=\frac{(2\mu_x \mu_y +c_1)(2\sigma_{xy} +c_2)}{({\mu}_x^2 +{\mu}_y^2 +c_1)({\sigma}_x^2 + {\sigma}_y^2 +c_2)}\,,
\end{align*}
here $\mu_x, \mu_y, {\sigma}_x^2, {\sigma}_y^2, \sigma_{xy} $ are the average, variance and covariance of $X$ and $Y$, respectively. The variables $c_1$ and $c_2$ are used to stabilize the division with weak denominator. To find overall performance of the image, we used a mean SSIM index (MSSIM),
\begin{align*}
\text{MSSIM}(X,Y)=\frac{1}{\text{N}} \sum_{i=1}^{\text{N}} \text{SSIM}(X_i,Y_i)\,,
\end{align*}
where $X$ and $Y$ are used to represent the original and the denoised images, respectively. $\text{N}$ is the number of local windows in the image, and $X_i$ and $Y_i$ are the image contents at the $i^{th}$ local window.\\
(d). Improvement in signal to noise ratio (ISNR) \cite{gonzalez2002digital}, is simply the difference between the improved and the original signal to noise ratio,
\begin{align*}
\text{ISNR} = 10\text{log}_{10}\left(\frac{ {\sum\limits_{i=1}^{\text{M}}} {\sum\limits_{j=1}^\text{N}} (I-I_0)^2  }{ \sum\limits_{i=1}^\text{M} \sum\limits_{j=1}^\text{N} (I-I_t)^2  }\right).
\end{align*}
To test the effectiveness of the proposed model, figure[\ref{fig:walkbridge_40}] presents the denoising results for Walkbridge test image that contain additive Gaussian noise of $\sigma=40$, which make the features hard to visualize. The restored outputs obtained from the state-of-the-art diffusion models depict that the denoised images are not satisfactory and the texture information of the image is ruined. From the figure[\ref{fig:3j}], it is easy to perceive that the proposed model can effectively eliminate all noise particles and preserve the original structure of the image.

\begin{figure}
        \centering 
        \begin{subfigure}[b]{0.3\textwidth}           
                \includegraphics[scale=0.2]{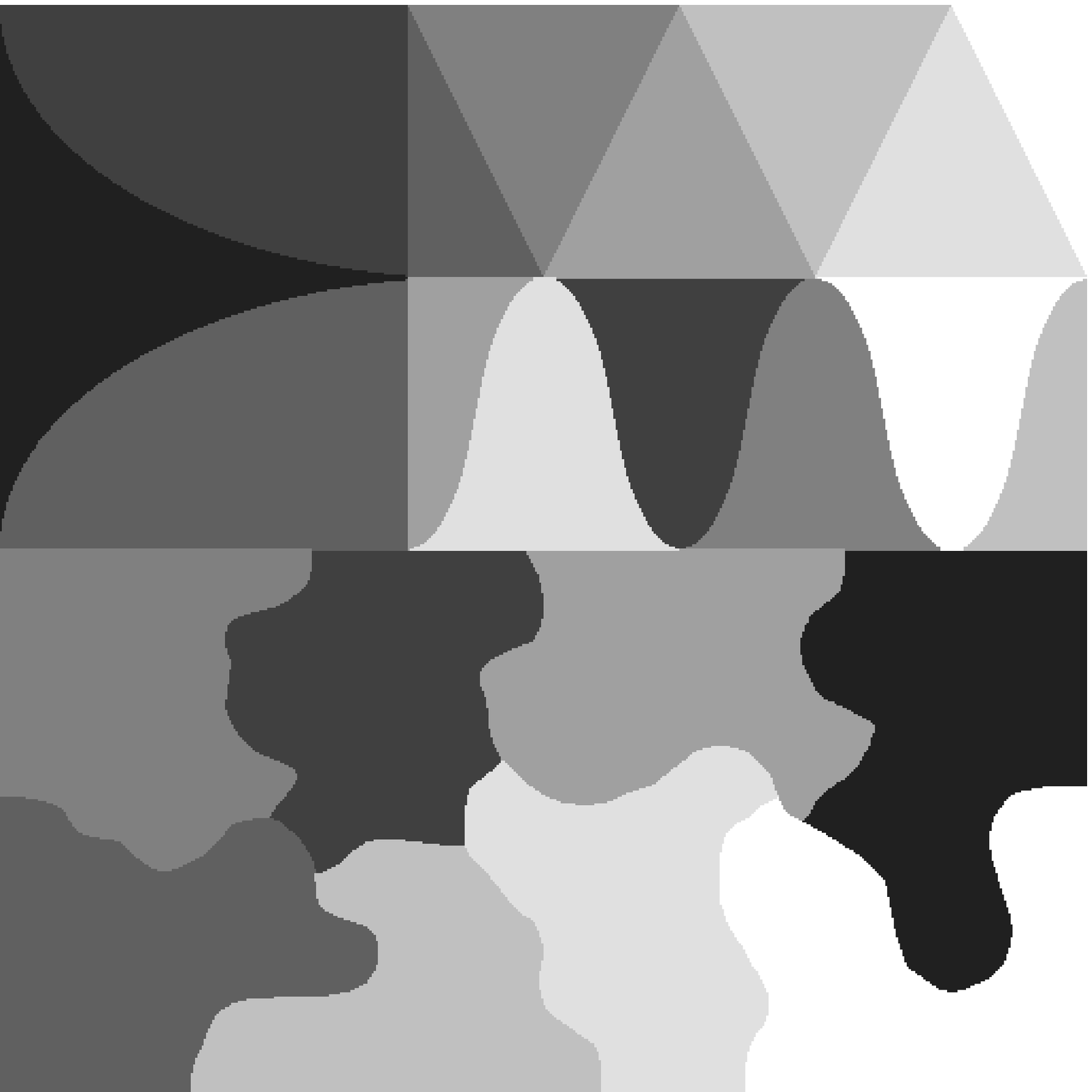}               
                \caption{Clear}
                \label{fig:6a}
        \end{subfigure}%
       \begin{subfigure}[b]{0.3\textwidth}
                \includegraphics[scale=0.2]{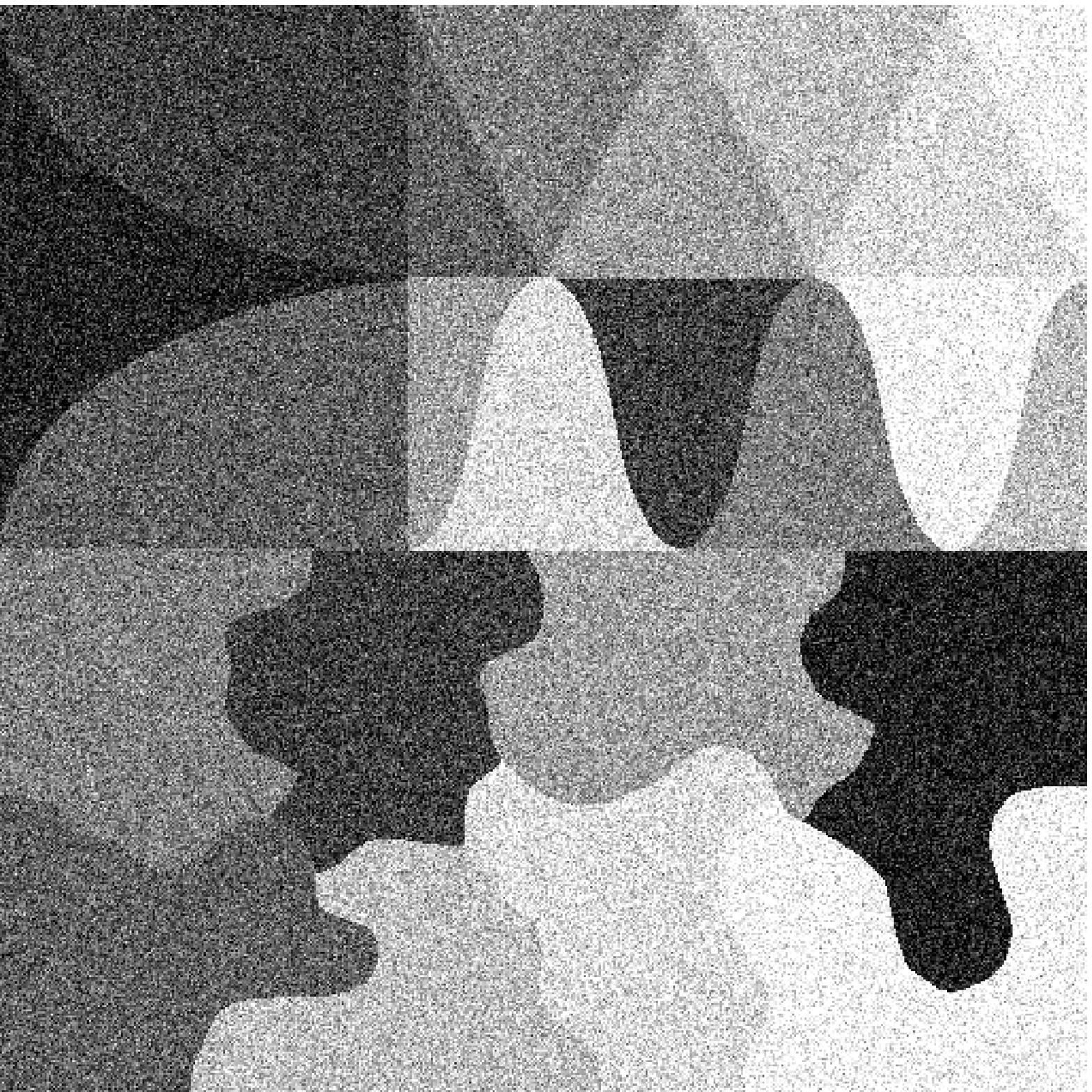}
                \caption{Noisy}
                \label{fig:6b}
        \end{subfigure}

         \begin{subfigure}[b]{0.3\textwidth}
                \includegraphics[scale=0.2]{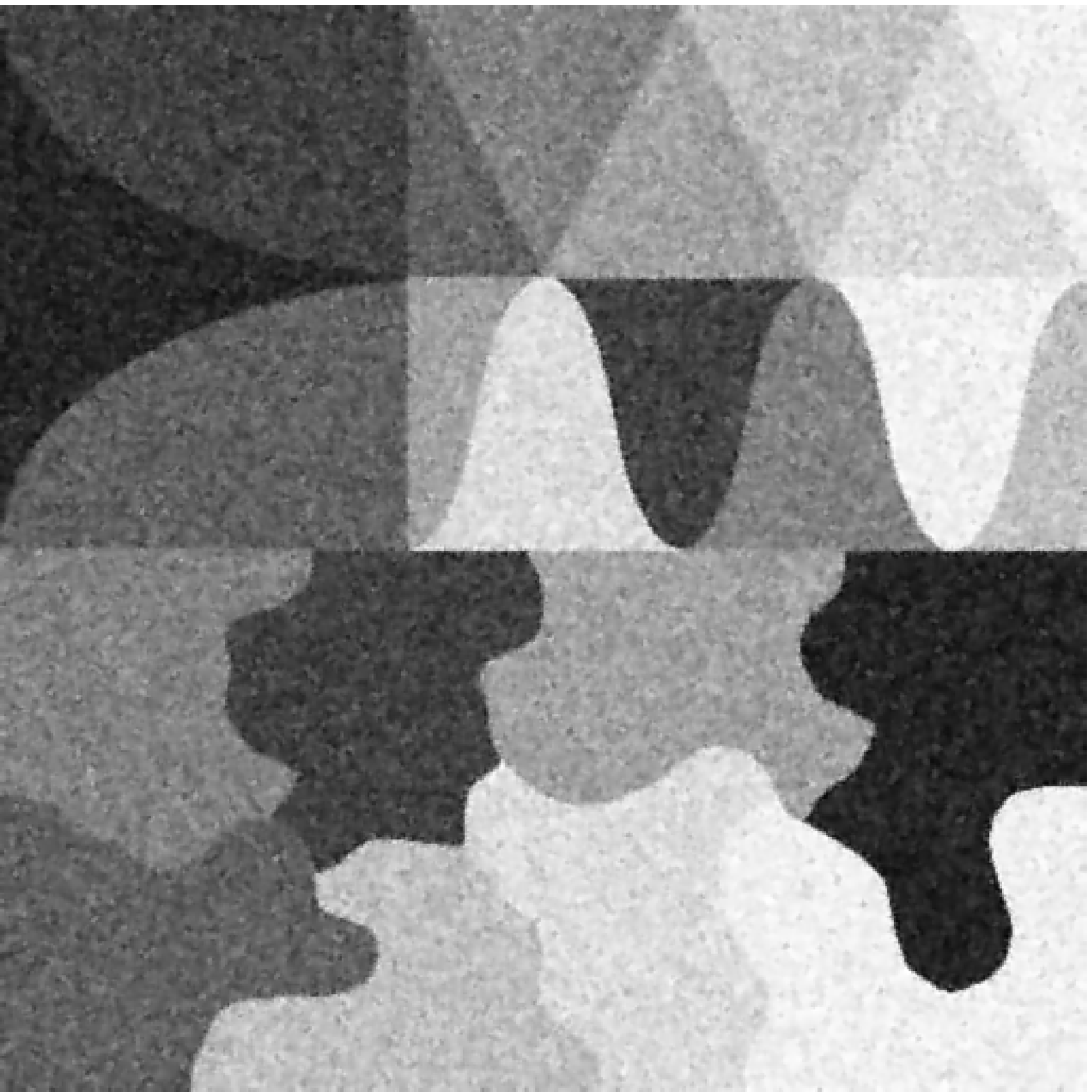}
                \caption{PV}
                \label{fig:6j}
        \end{subfigure}     
          \begin{subfigure}[b]{0.3\textwidth}
                \includegraphics[scale=0.2]{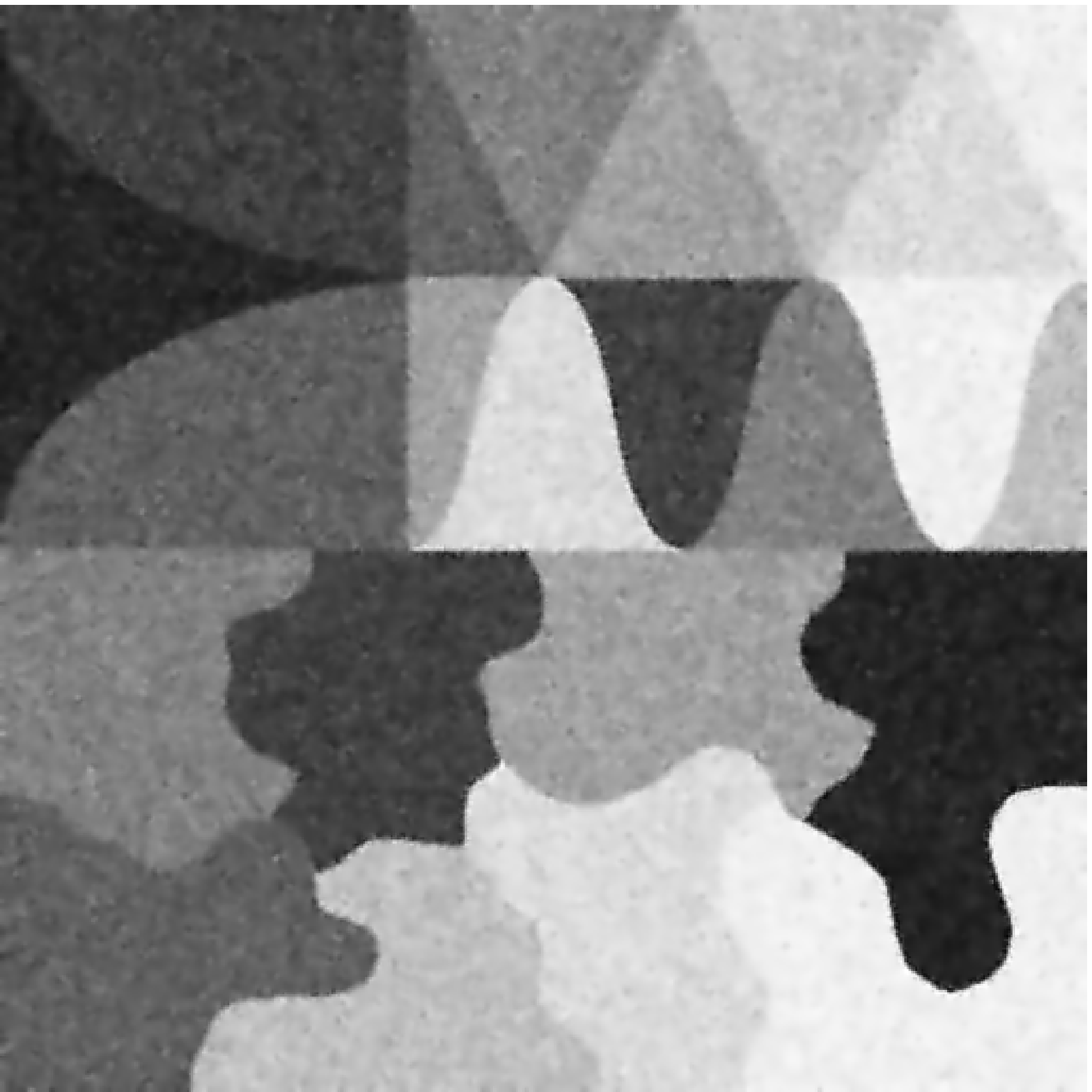}
                \caption{SYS}
                \label{fig:6m}
        \end{subfigure}
           \begin{subfigure}[b]{0.3\textwidth}
                \includegraphics[scale=0.2]{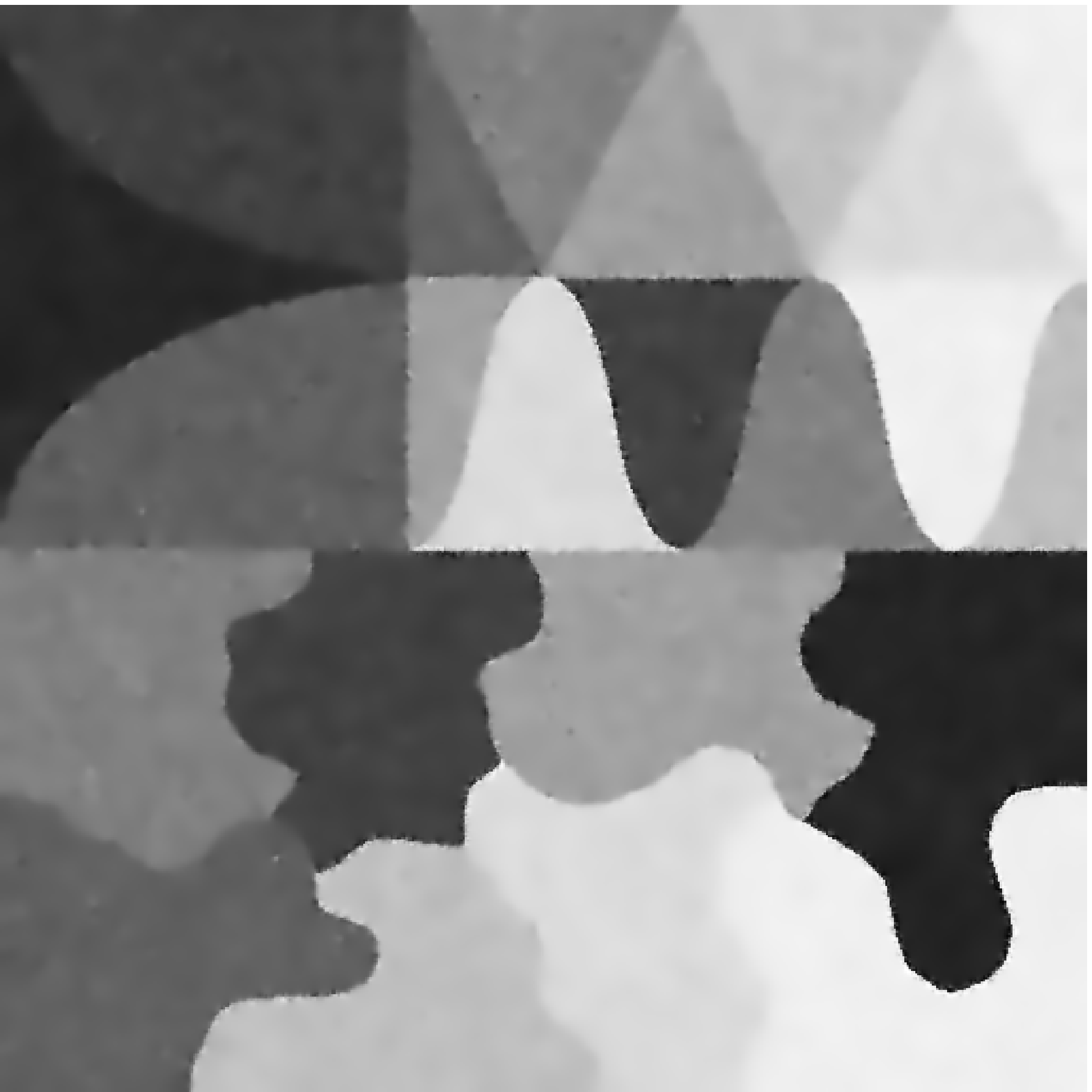}
                \caption{Proposed}
                \label{fig:6p}
        \end{subfigure}

        \begin{subfigure}[b]{0.3\textwidth}
                \includegraphics[scale=0.24]{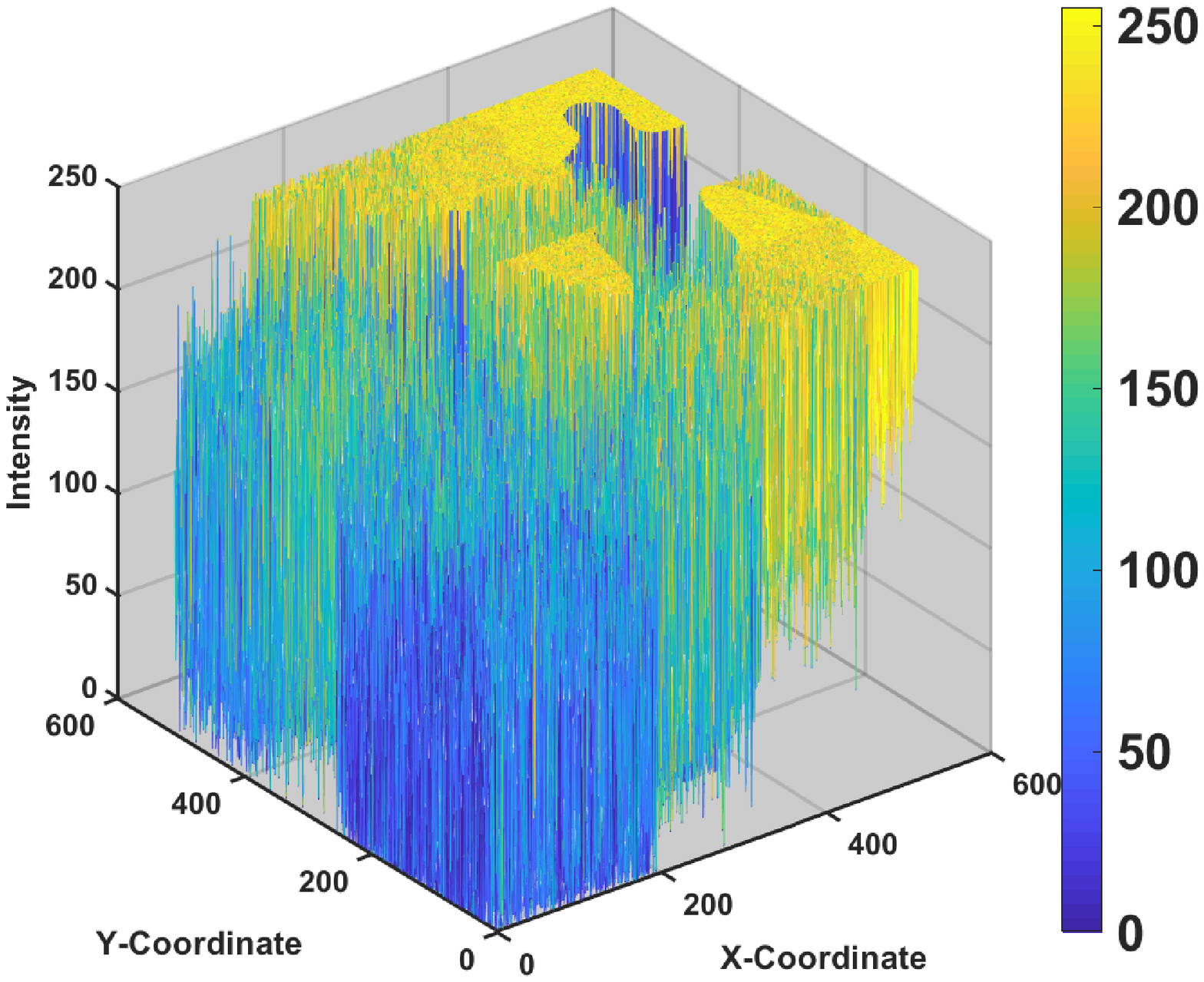}
                \caption{Noisy}
                \label{fig:6c}
        \end{subfigure} 
        
         \begin{subfigure}[b]{0.3\textwidth}
                \includegraphics[scale=0.24]{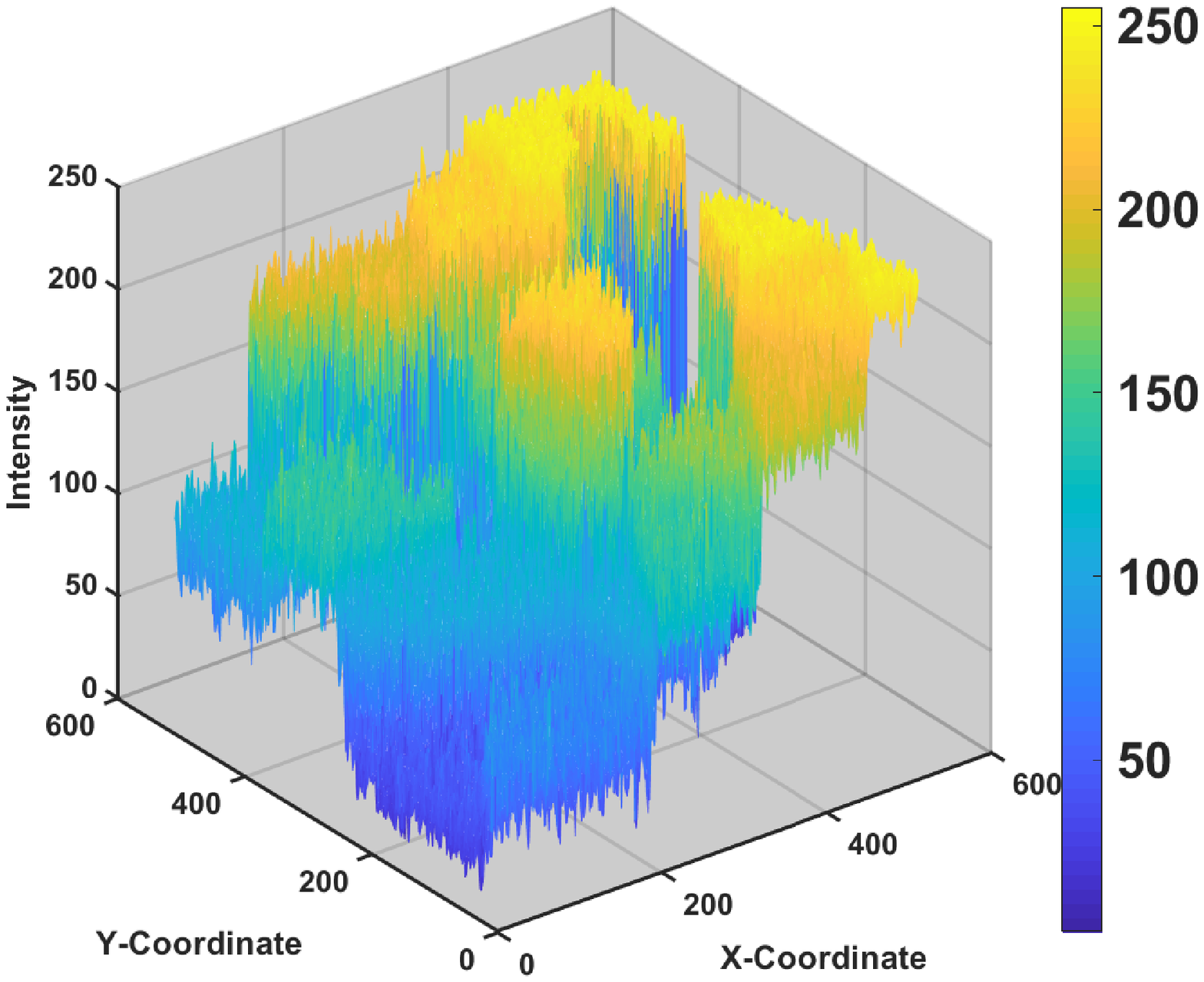}
                \caption{PV}
                \label{fig:6l}
        \end{subfigure}
           \begin{subfigure}[b]{0.3\textwidth}
                \includegraphics[scale=0.24]{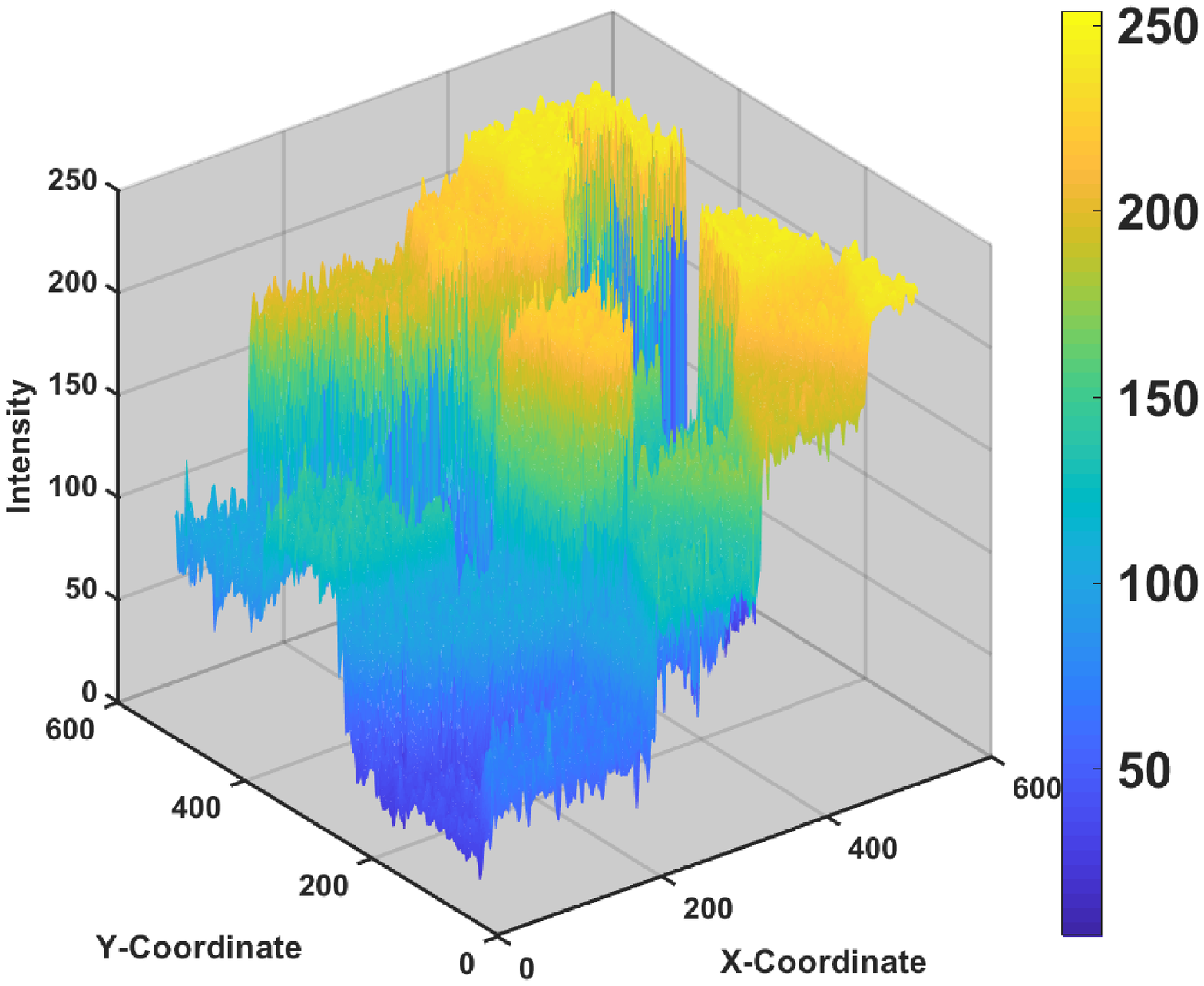}
                \caption{SYS}
                \label{fig:6o}
        \end{subfigure}
        \begin{subfigure}[b]{0.3\textwidth}
                \includegraphics[scale=0.24]{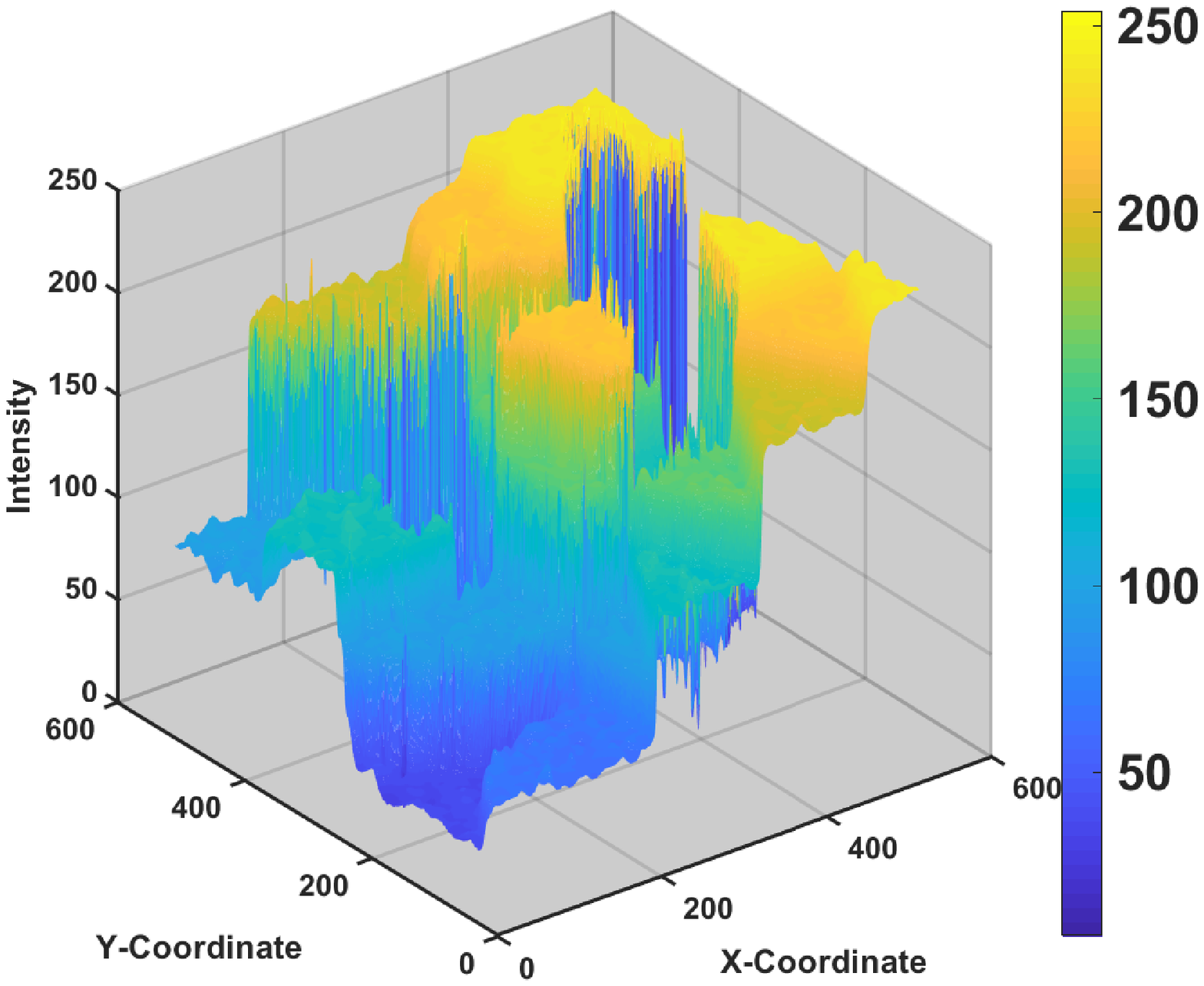}
                \caption{Proposed}
                \label{fig:6r}
        \end{subfigure}

\caption{(a) Clear image (b) Noisy image with Gaussian noise of mean 0.0 and s.d. 40 (c) PV Model, $\lambda=0.01, K=1$; (d) SYS Model, $\lambda=0.1, K=4$;  (e) Proposed Model, $\psi=1, k=3.75$.
 (f-i) 3D surfaces of images }\label{fig:texture_40_ratio}
\end{figure}

Further, to showcase the efficiency of our model, in figure \ref{fig:texture_40_ratio}, we compare the visual quality of the mosaic image computed by the present model along with the other existing models. Here we present the denoised images and 3D surface plots of the denoised images side by side for better comparison. From the figure, it can be observed that our model works better in terms of noise removal as compared to other models. From the 3D surface plots of the denoised images, one can see that our model leave fewer fluctuations with compare to the other models, it indicates that the proposed model not only removes noise efficiently but also preserves the fine structure as compared to the other models.

\begin{figure}
        \centering 
           \begin{subfigure}[b]{0.31\textwidth}           
                \includegraphics[scale=0.24]{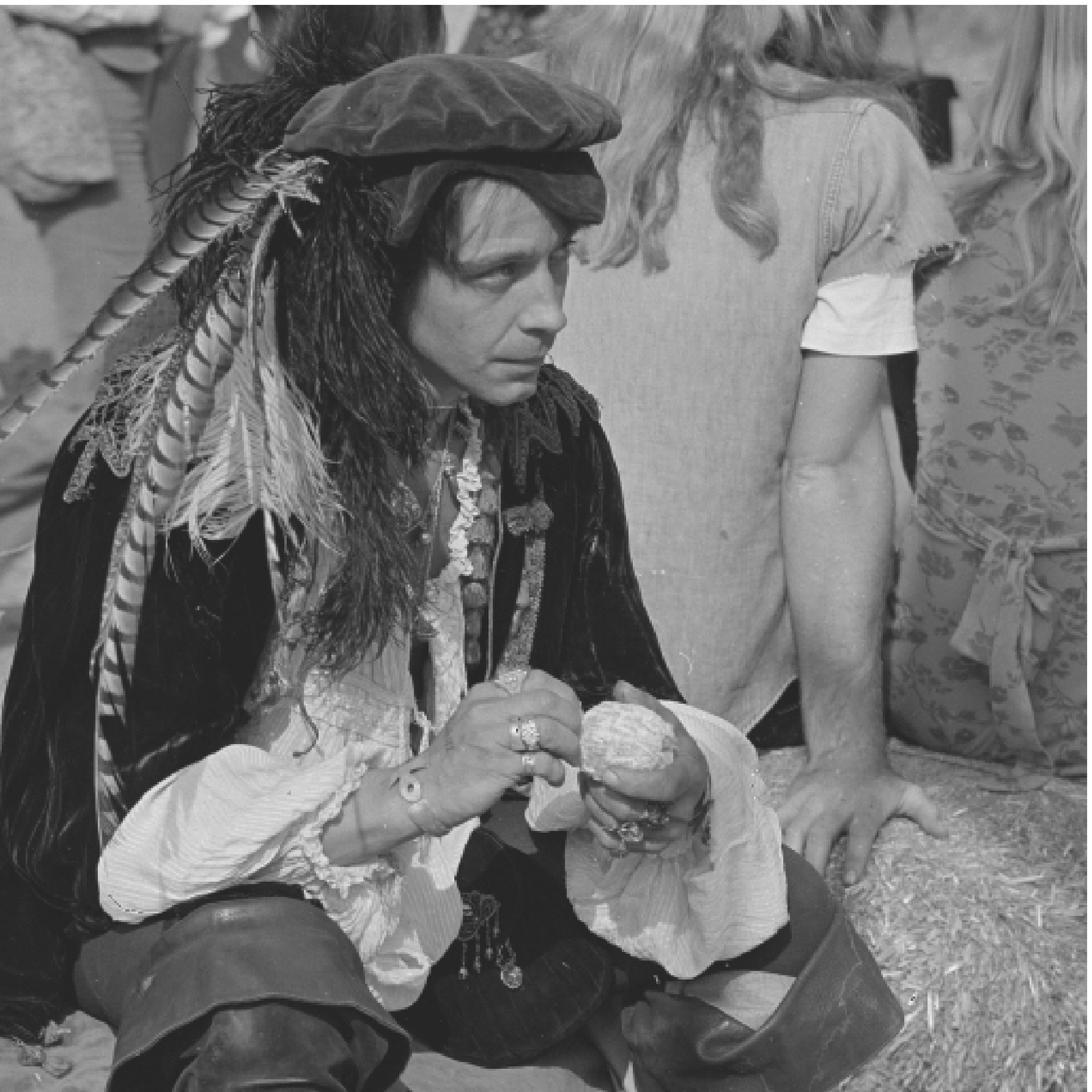}               
                \caption{Original}
                \label{fig:7a}
        \end{subfigure}%
        \begin{subfigure}[b]{0.31\textwidth}           
                \includegraphics[scale=0.24]{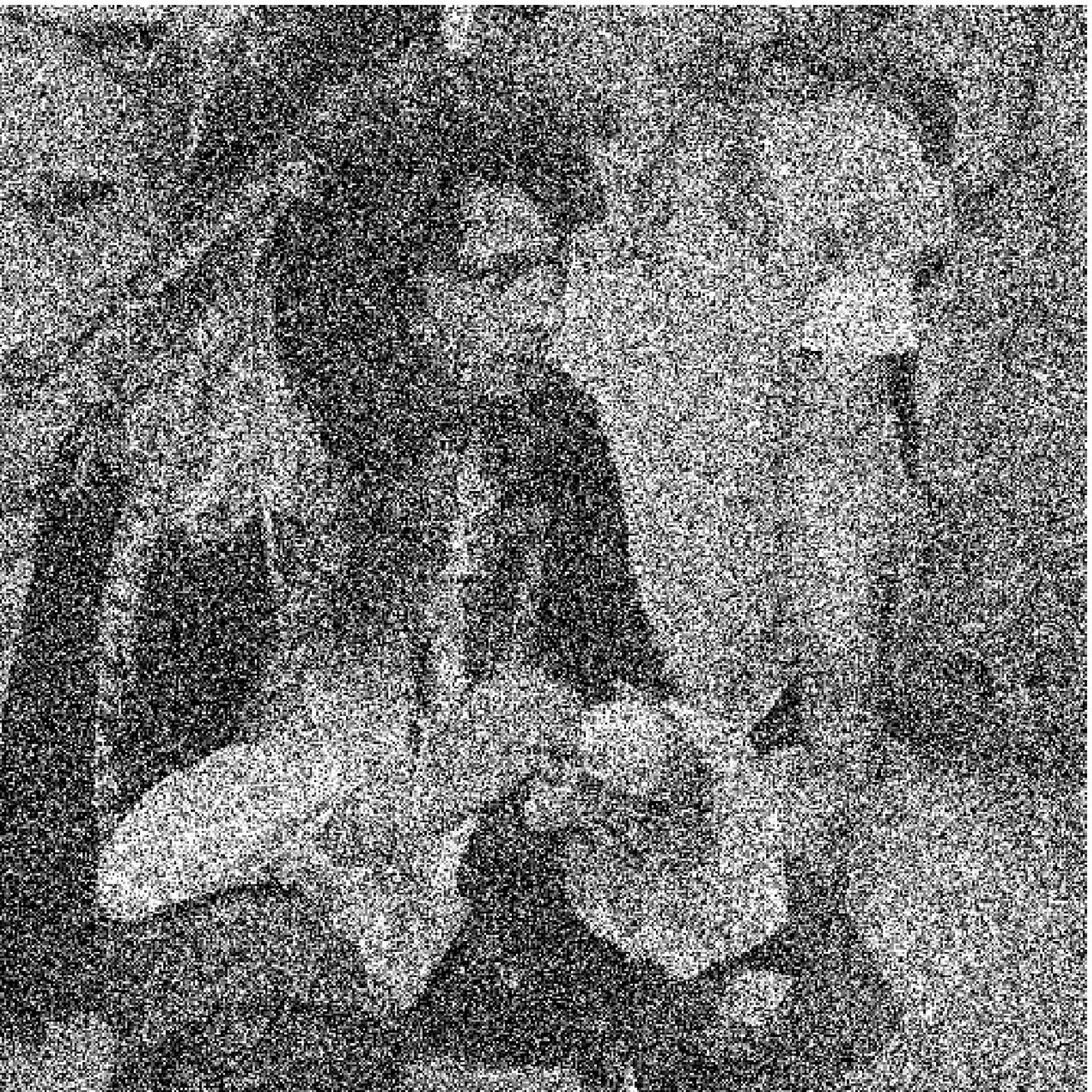}               
                \caption{Noisy}
                \label{fig:7b}
        \end{subfigure}%

       \begin{subfigure}[b]{0.3\textwidth}
                \includegraphics[scale=0.24]{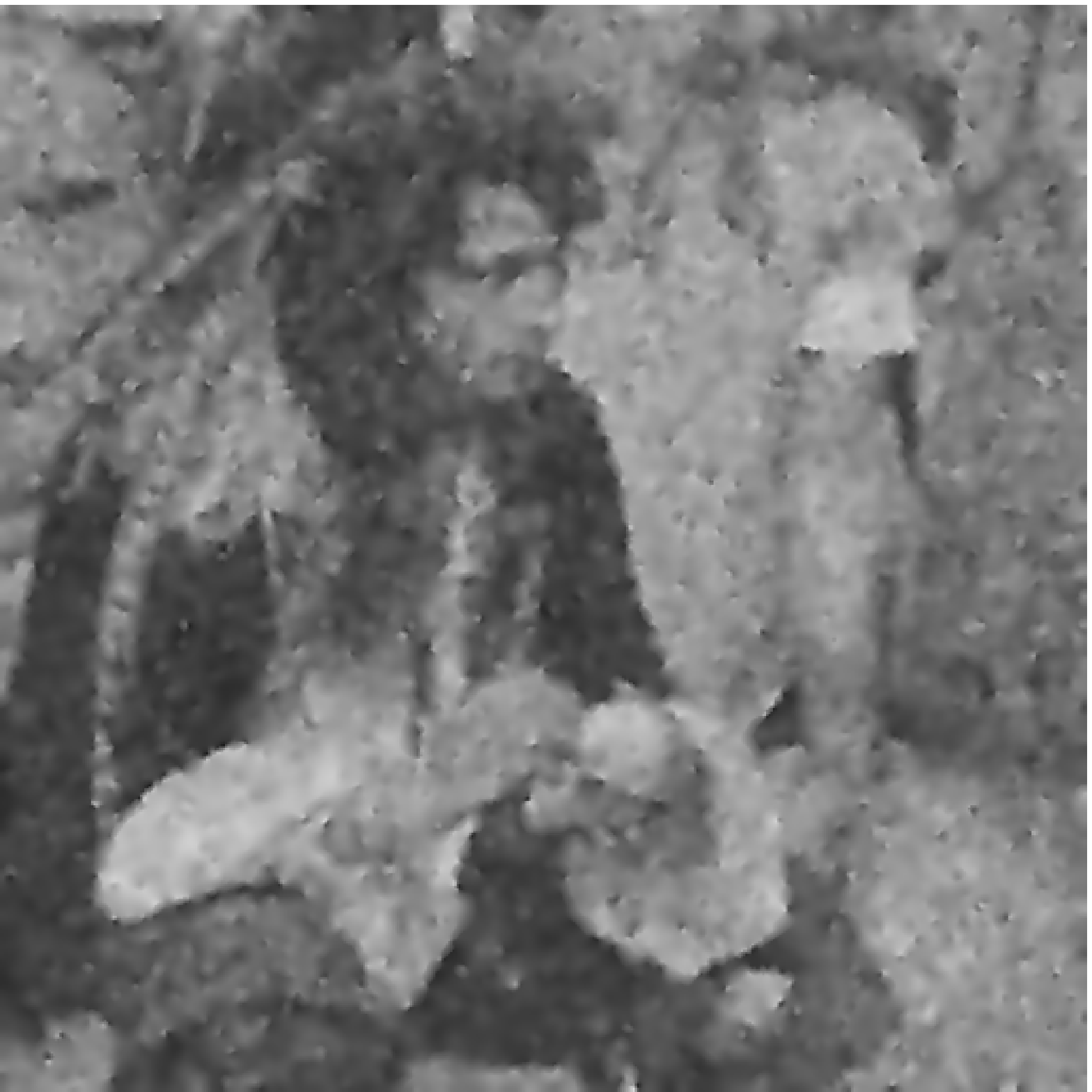}
                \caption{PV Model}
                \label{fig:7d}
        \end{subfigure}
        \begin{subfigure}[b]{0.3\textwidth}
                \includegraphics[scale=0.24]{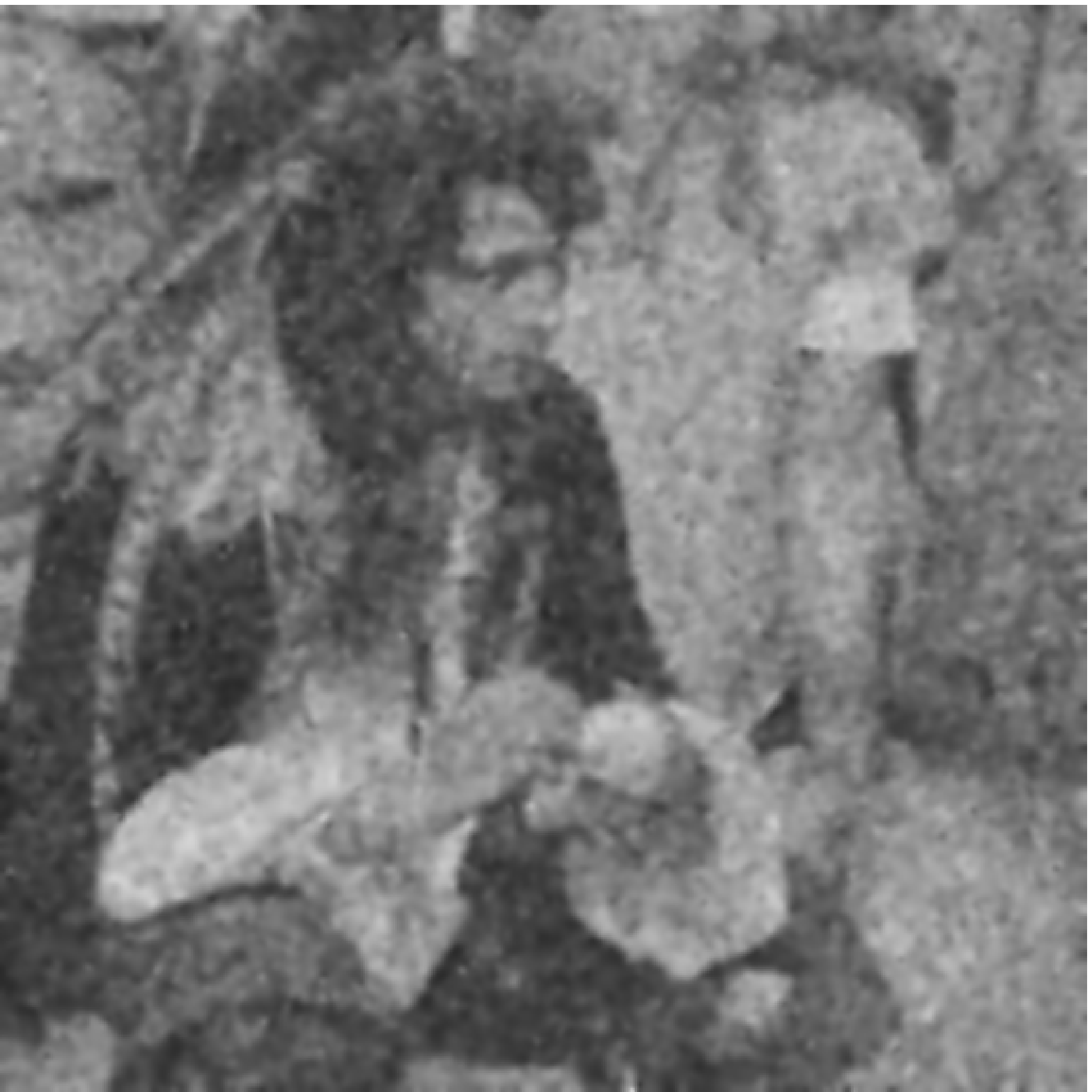}
                \caption{SYS Model}
                \label{fig:7e}
        \end{subfigure}      
        \begin{subfigure}[b]{0.3\textwidth}
                \includegraphics[scale=0.24]{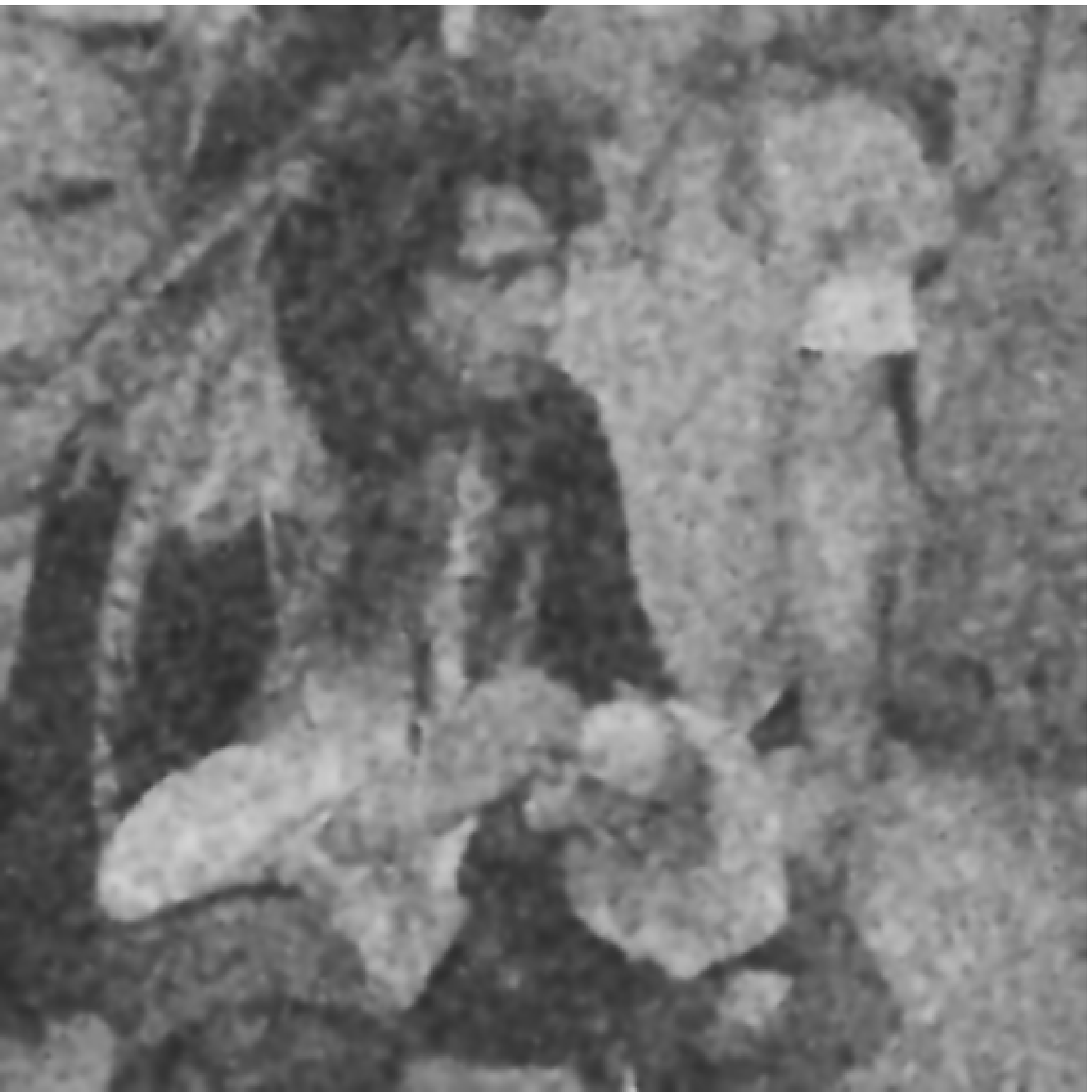}
                \caption{Proposed Model}
                \label{fig:7f}
        \end{subfigure}
\caption{ (a) Clear image (b)Heavily noised image with low SNR=3.70; Improved SNR values by different models: (c) SNR=15.41; $\lambda=0.01,K=11$ (d) SNR=15.37; $\lambda=1,K=10$ (e) SNR=15.53; $\psi=1,k=14$.}\label{heavy_noise_pirate_100}
\end{figure}

To further confirm the ability of the proposed model, in figure[\ref{heavy_noise_pirate_100}], we present the comparison of the denoised Pirate image which was initially degraded with additive Gaussian noise with very low SNR value $3.70$. From figure[\ref{heavy_noise_pirate_100}], it is clear that the proposed algorithm retains more texture information in addition to noise reduction, which makes our filtered image better in comparison to other models considered here.

\begin{figure}[]
        \centering
        \hspace{-1.5cm} 
        \begin{subfigure}[b]{0.28\textwidth}           
                \includegraphics[scale=0.22]{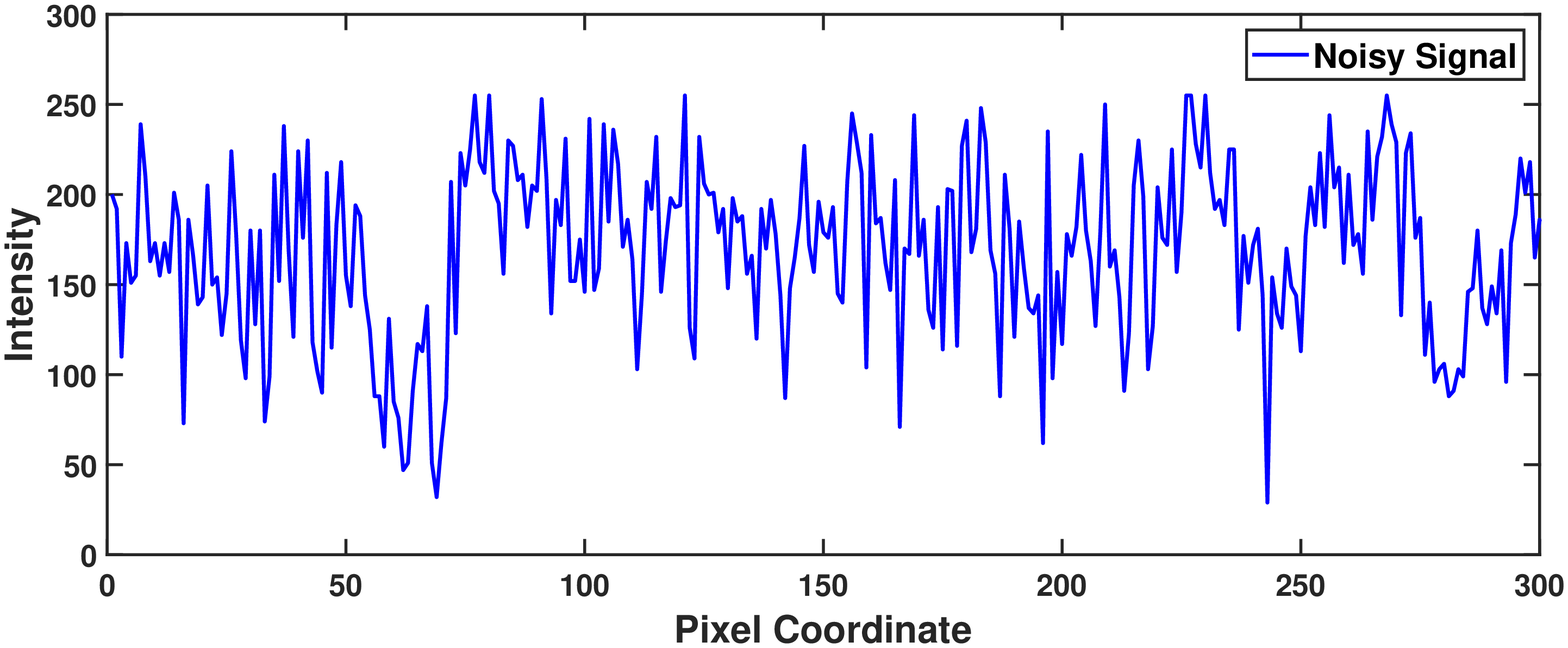}     
                \caption{}
                \label{fig:9a}
        \end{subfigure}%
        \hspace{2.5cm}
         \begin{subfigure}[b]{0.28\textwidth}           
                \includegraphics[scale=0.22]{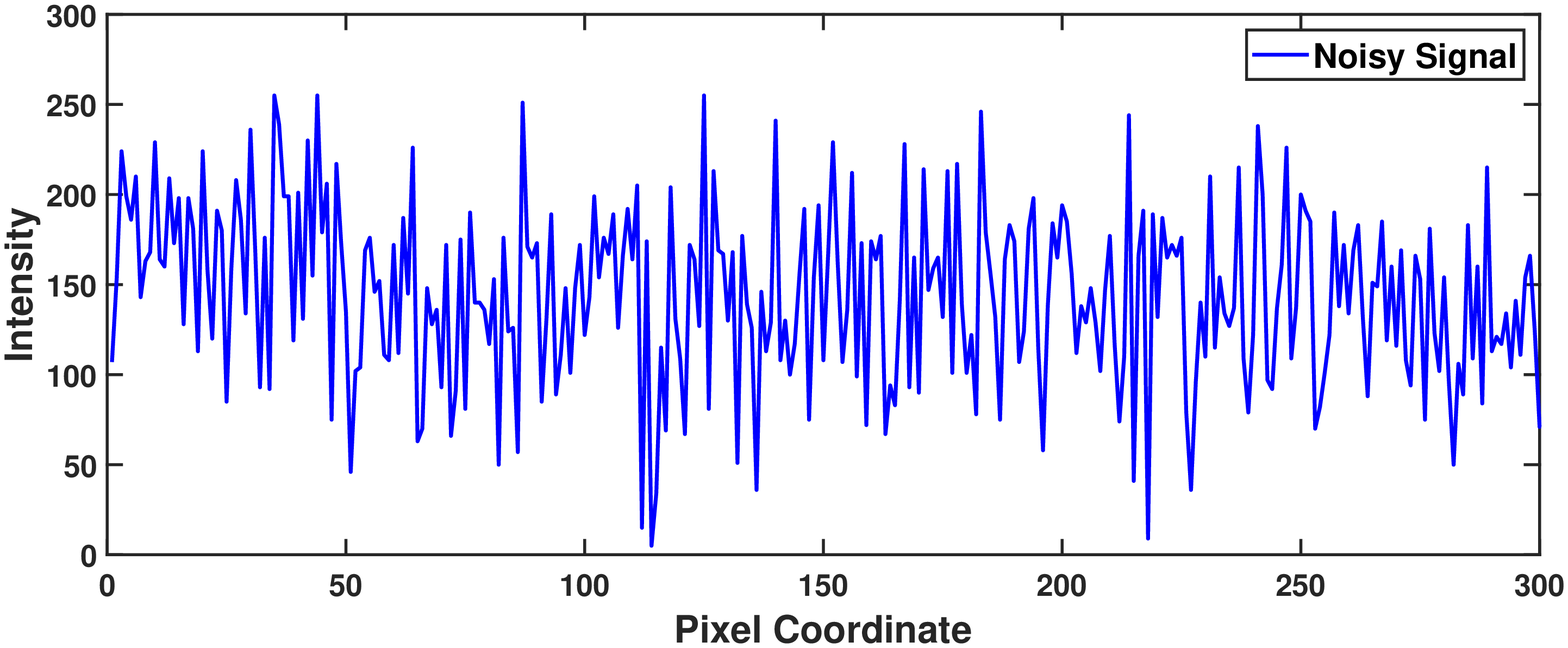}          
                \caption{}
                \label{fig:9b}
        \end{subfigure}%

         \hspace{-1.5cm}
         \begin{subfigure}[b]{0.28\textwidth}
                \includegraphics[scale=0.22]{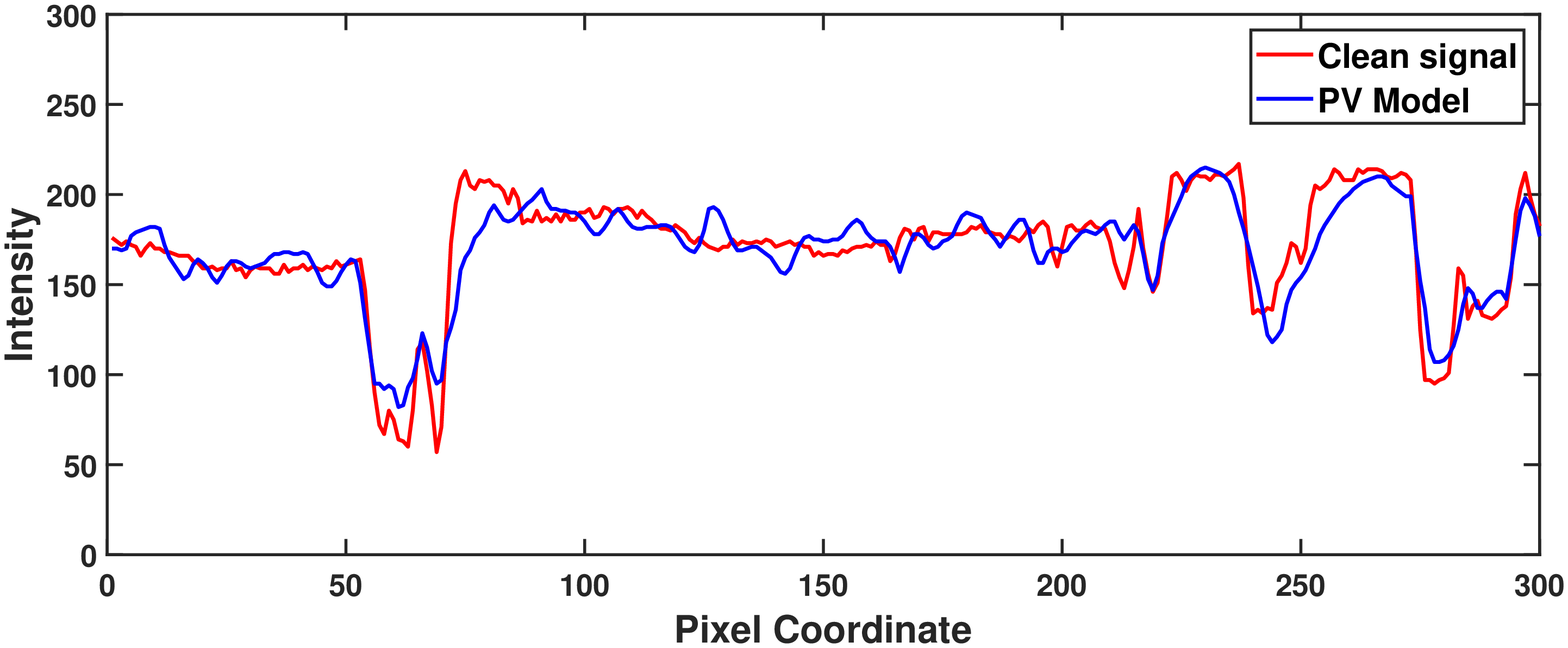}
                \caption{$\lambda=0.1, K=5$}
                \label{fig:9g}
        \end{subfigure}
        \hspace{2.5cm}
         \begin{subfigure}[b]{0.28\textwidth}           
                \includegraphics[scale=0.22]{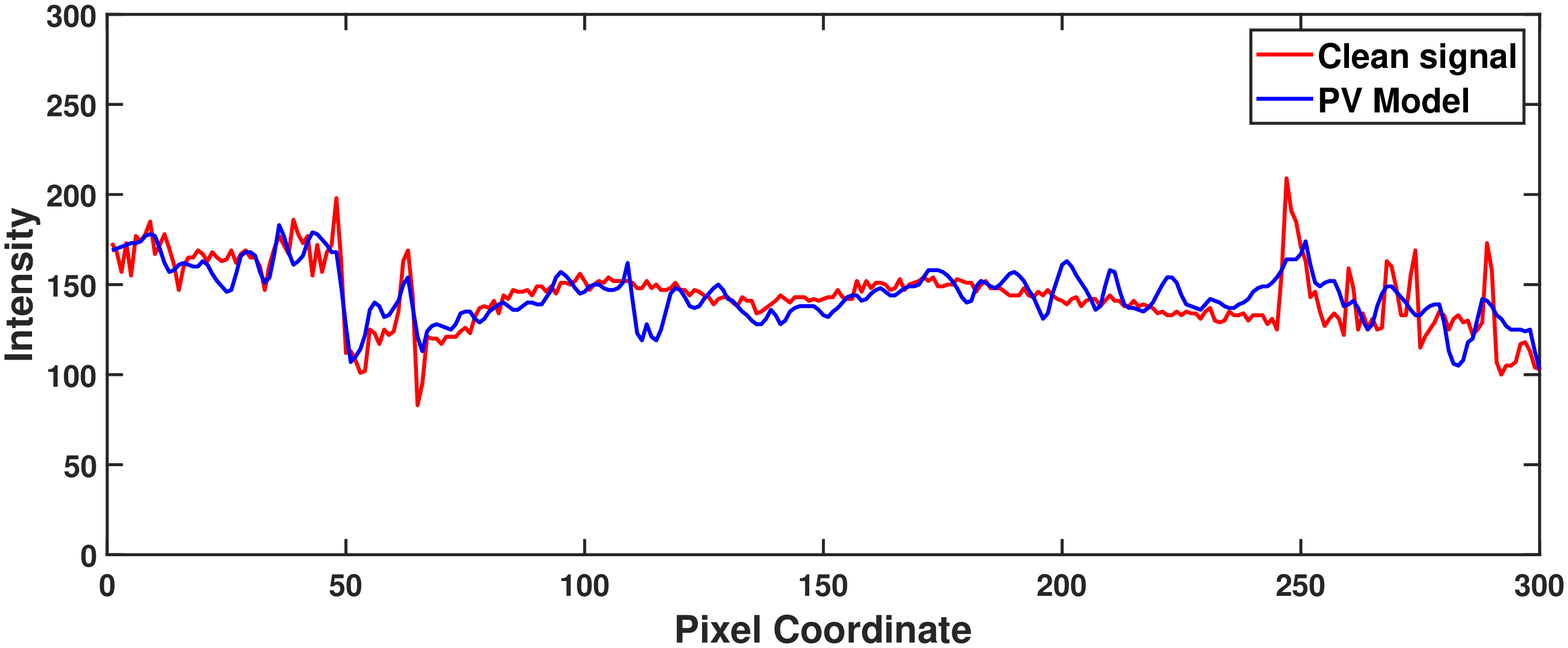}            
                \caption{$\lambda=0.1, K=15$}
                \label{fig:9h}
        \end{subfigure}
        
          \hspace{-1.5cm}
         \begin{subfigure}[b]{0.28\textwidth}
                \includegraphics[scale=0.22]{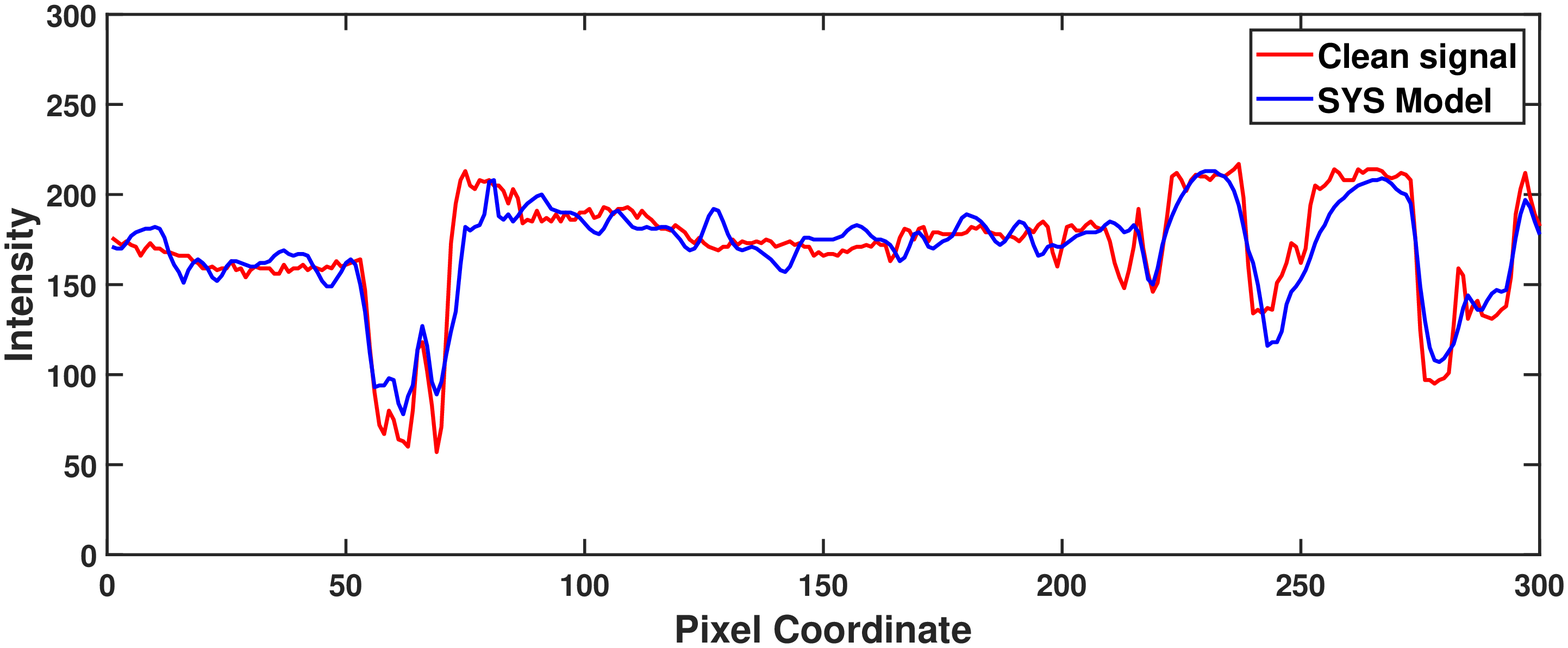}
                \caption{$\lambda=2, K=6$}
                \label{fig:9i}
        \end{subfigure}
        \hspace{2.5cm}
         \begin{subfigure}[b]{0.28\textwidth}           
                \includegraphics[scale=0.22]{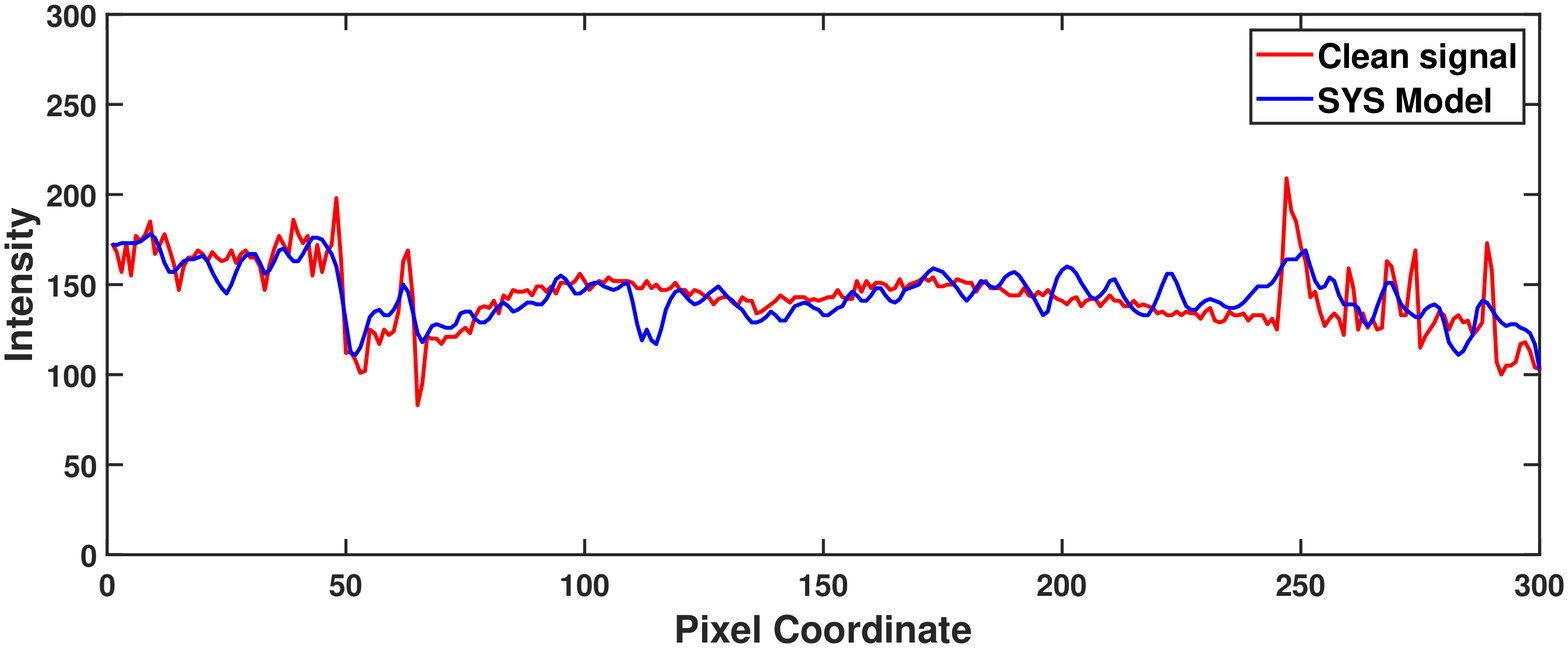}            
                \caption{$\lambda=4, K=9$}
                \label{fig:9j}
        \end{subfigure}

             \hspace{-1.5cm}
         \begin{subfigure}[b]{0.28\textwidth}
                \includegraphics[scale=0.22]{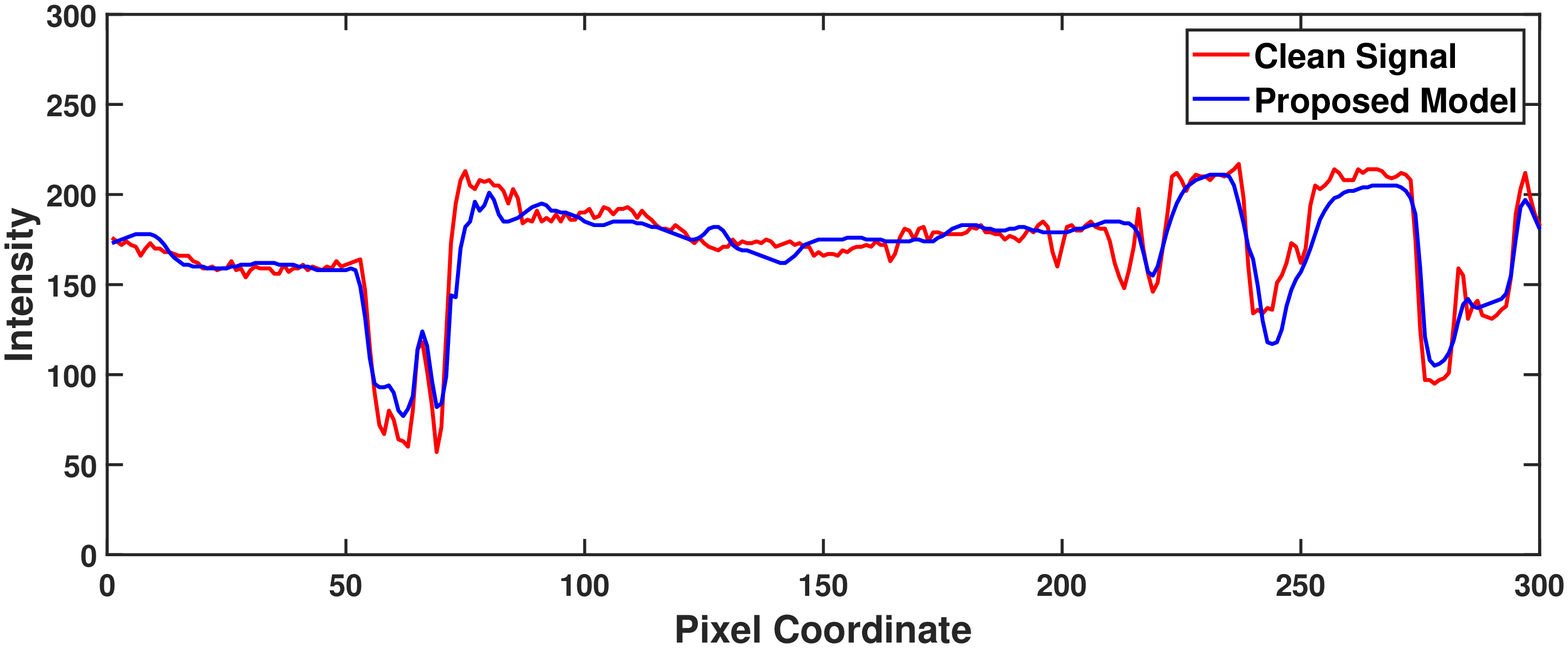}
                \caption{$\psi=1, k=4$}
                \label{fig:9k}
        \end{subfigure}
        \hspace{2.5cm}
         \begin{subfigure}[b]{0.28\textwidth}           
                \includegraphics[scale=0.22]{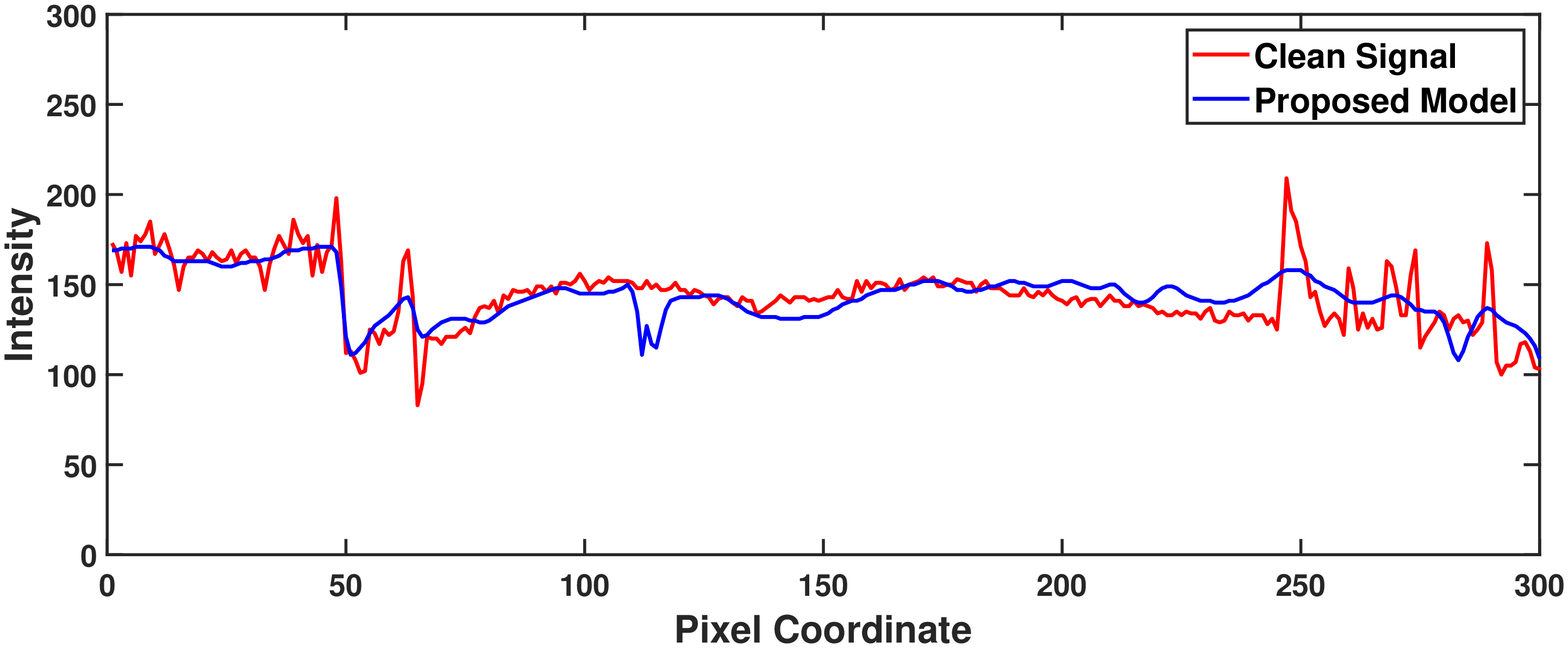}            
                \caption{$\psi=1, k=4.45$}
                \label{fig:9l}
        \end{subfigure}

        \caption{(a) The $219^{th}$ slice of Boat image corrupted with Gaussian noise of mean $0.0$ and $\sigma=40$; (b) The $95^{th}$ slice of Livingroom image corrupted with additive Gaussian noise of mean $0.0$ and $\sigma=50$; denoised using (c-d) PV model; (e-f) SYS model, and (g-h) Proposed model.}\label{fig:boat_40_livingroom_50_signal}
\end{figure}

Along with the qualitative analysis using full image surface, we also study the quality of the resultant images considering a slice of the Boat image and the Livingroom image with different noise levels. Figure \ref{fig:boat_40_livingroom_50_signal} shows the signals of the original, noisy and restored images obtained by the proposed CPDE model and other discussed CPDE models. From these figures, it is easy to conclude that the restored signals computed by the proposed model are more closure to the clean signals in comparison to other discussed models.


\begin{table}[]
\centering
\caption{Comparison of MSSIM and PSNR values for various approaches and proposed model. Clean image is degraded by additive Gaussian noise of mean $0.0$ and $\sigma=20$.}
\label{tab:table20}
\scalebox{0.6}{
        \begin{tabular}{cccccccccc}
\toprule
Images &Measure  &TV\cite{rudin1992nonlinear} &NS\cite{nitzberg1992nonlinear} &Luo\cite{luo2006coupled} &RD\cite{guo2011reaction} &NLM\cite{buades2005non}& PV\cite{prasath2014system}  & SYS\cite{sun2016class}  &Proposed\\
\midrule

Boat        &MSSIM        &$0.8232$     &$0.8883$       &$0.8774$ &$0.8886$ & $0.9105$     & $0.8900 $  & $0.8943$ &    $\mathbf{0.9219}$ \\
            &PSNR         &$28.10$      &$28.81$        &$29.40$  &$29.61$  & $30.21$      & $29.69  $  & $29.84 $ &    $\mathbf{30.24}$ \\
Lake        &MSSIM        &$0.8539$     &$0.8996$       &$0.8775$ &$0.8788$ & $0.9140$     & $0.9060 $  & $0.9074$ &    $\mathbf{0.9268}$ \\
            &PSNR         &$28.37$      &$28.07$        &$28.65$  &$28.61$  & $29.27$      & $28.87  $  & $29.08$  &    $\mathbf{29.57}$ \\
Livingroom  &MSSIM        &$0.8790$     &$0.8698$       &$0.8855$ &$0.8849$ & $0.8870$     & $0.8903 $  & $0.8912$ &    $\mathbf{0.9019}$ \\
            &PSNR         &$28.51$      &$28.62$        &$28.68$  &$28.80$  & $29.11$      & $28.87  $  & $28.93$  &    $\mathbf{29.22}$ \\
Mandril     &MSSIM        &$0.8958$     &$0.9001$       &$0.8909$ &$0.8997$ & $0.8684$     & $0.9015 $  & $0.9029$ &    $\mathbf{0.9034}$ \\
            &PSNR         &$27.51$      &$27.71$        &$27.77$  &$27.95$  & $27.62$      & $28.24  $  & $28.37$  &    $\mathbf{28.46}$ \\
Pirate      &MSSIM        &$0.8615$     &$0.9053$       &$0.8764$ &$0.8714$ & $0.8908$     & $0.8936 $  & $0.8943$ &    $\mathbf{0.9055}$ \\
            &PSNR         &$28.54$      &$29.40$        &$29.20$  &$29.07$  & $29.69$      & $29.51  $  & $29.53$  &    $\mathbf{29.90}$ \\
Walkbridge  &MSSI         &$0.9077$     &$0.9141$       &$0.9070$ &$0.9101$ & $0.8747$     & $0.9093 $  & $0.9089$ &    $\mathbf{0.9143}$ \\
            &PSNR         &$26.79$      &$26.83$        &$27.05$  &$26.99$  & $26.72$      & $26.99  $  & $27.02$  &    $\mathbf{27.24}$ \\
Woman       &MSSIM        &$0.8007$     &$0.9191$       &$0.8257$ &$0.8206$ & $0.9243$     & $0.9262 $  & $0.9225$ &    $\mathbf{0.9396}$ \\
            &PSNR         &$30.35$      &$33.83$        &$31.31$  &$31.29$  & $33.64$      & $34.47  $  & $34.59 $ &    $\mathbf{35.03}$ \\
Mosaic      &MSSIM        &$0.9699$     &$0.9091$       &$0.9270$ &$0.9600$ & $0.9588$     & $0.9675 $  & $0.9290$ &    $\mathbf{0.9780}$ \\
            &PSNR         &$34.28$      &$32.71$        &$33.66$  &$34.44$  & $34.12$      & $33.38  $  & $32.54$  &    $\mathbf{34.74}$ \\               
\bottomrule
\end{tabular}
}

\end{table}
\begin{table}[]
\centering
  \vspace{1.0cm}
        \caption{Comparison of MSSIM and PSNR values for various approaches and proposed model. Clean image is degraded by additive Gaussian noise of mean $0.0$ and $\sigma=30$.}
\label{tab:table30}
       \scalebox{0.6}{
        \begin{tabular}{cccccccccc}
\toprule
Images &Measure  &TV\cite{rudin1992nonlinear} &NS\cite{nitzberg1992nonlinear} &Luo\cite{luo2006coupled} &RD\cite{guo2011reaction} &NLM\cite{buades2005non} & PV\cite{prasath2014system}  & SYS\cite{sun2016class}  &Proposed\\
\midrule

Boat         &MSSIM          & $0.7157$ & $0.6923$ & $0.7324$ & $0.7327$ & $0.8568$  & $0.8360$ & $0.8400$ & $\mathbf{0.8841}$ \\
             &PSNR           & $25.28$  & $25.71$  & $25.83$  & $26.01$  & $28.02$   & $27.79$  & $27.96 $ & $\mathbf{28.47}$  \\
Lake         &MSSIM          & $0.7644$ & $0.7676$ & $0.7751$ & $0.7820$ & $0.8705$  & $0.8576$ & $0.8637$ & $\mathbf{0.8956}$ \\
             &PSNR           & $25.27$  & $25.41$  & $25.72$  & $26.01$  & $27.38$   & $27.12 $ & $27.33$  & $\mathbf{27.87}$ \\
Livingroom   &MSSIM          & $0.7859$ & $0.8510$ & $0.7922$ & $0.7969$ & $0.8229$  & $0.8273$ & $0.8321$ & $\mathbf{0.8515}$ \\
             &PSNR           & $25.91$  & $27.28$  & $26.14$  & $26.40$  & $26.90$   & $26.97 $ & $27.02$  & $\mathbf{27.39}$ \\
Mandril      &MSSIM          & $0.8188$ & $0.7633$ & $0.8262$ & $0.8268$ & $0.7873$  & $0.8309$ & $0.8395$ & $\mathbf{0.8454}$ \\
             &PSNR           & $24.94$  & $25.46$  & $25.55$  & $25.52$  & $25.43$   & $26.13 $ & $26.27$  & $\mathbf{26.47}$ \\
Pirate       &MSSIM          & $0.7745$ & $0.8348$ & $0.7765$ & $0.7905$ & $0.8346$  & $0.8392$ & $0.8400$ & $\mathbf{0.8575}$ \\
             &PSNR           & $26.18$  & $27.05$  & $26.26$  & $26.86$  & $27.72$   & $27.86 $ & $27.81$  & $\mathbf{28.20}$ \\
Walkbridge   &MSSIM          & $0.8391$ & $0.8194$ & $0.8431$ & $0.8453$ & $0.8015$  & $0.8462$ & $0.8508$ & $\mathbf{0.8566}$ \\
             &PSNR           & $24.60$  & $24.07$  & $24.91$  & $24.99$  & $24.96$   & $25.19$  & $25.22$  & $\mathbf{25.47}$ \\
Woman        &MSSIM          & $0.6834$ & $0.8559$ & $0.6805$ & $0.7108$ & $0.8861$  & $0.8963$ & $0.8922$ & $\mathbf{0.9200}$ \\
             &PSNR           & $27.58$  & $31.16$  & $27.47$  & $28.58$  & $31.38$   & $32.58$  & $32.64 $ & $\mathbf{33.25}$ \\
Mosaic       &MSSIM          & $0.9558$ & $0.8827$ & $0.8800$ & $0.9467$ & $0.9477$  & $0.9520$ & $0.8764$ & $\mathbf{0.9659}$ \\
             &PSNR           & $31.63$  & $30.38$  & $30.68$  & $31.65$  & $31.69$   & $31.10$  & $ 30.24$ & $\mathbf{31.97}$ \\           
\bottomrule
\end{tabular}
}
\end{table}
\begin{table}[]
\centering
\caption{Comparison of MSSIM and PSNR values for various approaches and proposed model. Clean image is degraded by additive Gaussian noise of mean $0.0$ and $\sigma=40$.}
\label{tab:table40}
        \scalebox{0.6}{
        \begin{tabular}{cccccccccc}
\toprule
Images &Measure  &TV\cite{rudin1992nonlinear} &NS\cite{nitzberg1992nonlinear} &Luo\cite{luo2006coupled} &RD\cite{guo2011reaction} &NLM\cite{buades2005non} & PV\cite{prasath2014system}  & SYS\cite{sun2016class}  &Proposed\\
\midrule
Boat    &MSSIM          & $0.6305$ &$0.6146$ &$0.6290$  &$0.6481$  &$0.8022$    & $0.7862$ & $0.7798$  & $\mathbf{0.8470}$ \\
        &PSNR           & $23.43$  &$23.41$  &$23.35$   &$24.19$   &$26.40$     & $26.51$  & $26.60 $  & $\mathbf{27.21}$ \\
Lake    &MSSIM          & $0.6964$ &$0.7052$ &$0.6927$  &$0.7117$  &$0.8245$    & $0.8154$ & $0.8173$  & $\mathbf{0.8644}$ \\
        &PSNR           & $23.67$  &$23.70$  &$23.62$   &$24.33$   &$25.93$     & $25.94$  & $26.05$   & $\mathbf{26.62}$ \\
Livingroom &MSSIM       & $0.6941$ &$0.6989$ &$0.6855$  &$0.7095$  &$0.7600$    & $0.7699$ & $0.7731$  & $\mathbf{0.7991}$ \\
       &PSNR            & $23.85$  &$23.36$  &$23.56$   &$24.50$   &$25.37$     & $25.68$  & $25.70$   & $\mathbf{26.11}$ \\
Mandril &MSSIM          & $0.7409$ &$0.7456$ &$0.7462$  &$0.7525$  &$0.7164$    & $0.7692$ & $0.7777$  & $\mathbf{0.7844}$ \\
        &PSNR           & $23.14$  &$23.24$  &$23.41$   &$23.82$   &$24.01$     & $24.75$  & $24.90$   & $\mathbf{25.10}$ \\
Pirate  &MSSIM          & $0.6882$ &$0.6913$ &$0.6800$  &$0.7067$  &$0.7849$    & $0.7869$ & $0.7870$  & $\mathbf{0.8132}$ \\
        &PSNR           & $24.32$  &$24.23$  &$24.04$   &$25.06$   &$26.35$     & $26.64$  & $26.58$   & $\mathbf{27.06}$ \\
Walkbridge &MSSIM       & $0.7728$ &$0.7883$ &$0.7722$  &$0.7821$  &$0.7371$    & $0.7922$ & $0.7944$  & $\mathbf{0.7985}$ \\
        &PSNR           & $23.05$  &$23.48$  &$23.08$   &$23.54$   &$23.76$     & $24.01$  & $24.05$   & $\mathbf{24.24}$ \\
Woman   &MSSIM          & $0.5976$ &$0.7317$ &$0.5671$  &$0.6235$  &$0.8513$    & $0.8690$ & $0.8641$  & $\mathbf{0.9037}$ \\
        &PSNR           & $25.75$  &$27.59$  &$25.93$   &$26.61$   &$29.75$     & $31.00$  & $30.96$   & $\mathbf{31.65}$ \\
Mosaic  &MSSIM          & $0.9384$ &$0.8892$ &$0.8315$  &$0.9328$  &$0.9388$    & $0.9303$ & $0.8373$  & $\mathbf{0.9508}$ \\
         &PSNR          & $29.43$  &$28.51$  &$28.38$   &$29.46$   &$29.66$     & $29.04$  & $28.15$   & $\mathbf{29.77}$ \\          
        
\bottomrule
\end{tabular}
}
\end{table}
\begin{table}[]
\centering
\caption{Comparison of MSSIM and PSNR values for various approaches and proposed model. Clean image is degraded by additive Gaussian noise of mean $0.0$ and $\sigma=50$.}
\label{tab:table50}
        \scalebox{0.6}{
        \begin{tabular}{cccccccccc}
\toprule
Images &Measure  &TV \cite{rudin1992nonlinear} &NS \cite{nitzberg1992nonlinear} &Luo \cite{luo2006coupled} &RD \cite{guo2011reaction} &NLM\cite{buades2005non} & PV\cite{prasath2014system}  & SYS\cite{sun2016class}  &Proposed\\
\midrule
        
Boat        &MSSIM    & $0.5648$ & $0.5793$ &$0.5519$ &$0.5813$ &$0.7523$     & $0.7425$ &$ 0.7333$ & $\mathbf{0.8159}$ \\
            &PSNR     & $22.02$  & $21.70$  &$21.50$  &$22.75$  &$25.14$      & $25.48 $ &$ 25.46 $ & $\mathbf{26.16}$ \\
Lake        &MSSIM    & $0.6404$ & $0.6657$ &$0.6232$ &$0.6540$ &$0.7828$     & $0.7747$ &$ 0.7710$ & $\mathbf{0.8334}$ \\
            &PSNR     & $22.47$  & $22.09$  &$21.98$  &$23.06$  &$24.75$      & $24.89 $ &$ 24.94$  & $\mathbf{25.53}$ \\
Livingroom  &MSSIM    & $0.6265$ & $0.6424$ &$0.5981$ &$0.6429$ &$0.7090$     & $0.7245$ &$ 0.7229$ & $\mathbf{0.7562}$ \\
            &PSNR     & $22.51$  & $21.57$  &$21.54$  &$23.17$  &$24.28$      & $24.75 $ &$ 24.72 $ & $\mathbf{25.17}$ \\
Mandril     &MSSIM    & $0.6707$ & $0.6177$ &$0.6765$ &$0.6846$ &$0.6505$     & $0.7102$ &$ 0.7173$ & $\mathbf{0.7261}$ \\
            &PSNR     & $21.77$  & $21.48$  &$22.06$  &$22.49$  &$22.99$      & $23.69 $ &$ 23.87 $ & $\mathbf{24.06}$ \\
Pirate      &MSSIM    & $0.6275$ & $0.6486$ &$0.5969$ &$0.6450$ &$0.7423$     & $0.7483$ &$ 0.7454$ & $\mathbf{0.7762}$ \\
            &PSNR     & $23.11$  & $21.66$  &$22.15$  &$23.78$  &$25.24$      & $25.68 $ &$ 25.60 $ & $\mathbf{26.10}$ \\
Walkbridge  &MSSIM    & $0.7150$ & $0.7046$ &$0.7025$ &$0.7243$ &$0.6838$     & $0.7414$ &$ 0.7428$ & $\mathbf{0.7483}$ \\
            &PSNR     & $22.00$  & $21.75$  &$21.52$  &$22.44$  &$22.84$      & $23.10 $ &$ 23.16 $ & $\mathbf{23.34}$ \\
Woman       &MSSIM    & $0.5318$ & $0.5702$ &$0.5626$ &$0.6143$ &$0.8150$     & $0.8430$ &$ 0.8354$ & $\mathbf{0.8830}$ \\
            &PSNR     & $24.30$  & $22.02$  &$22.56$  &$25.02$  &$28.23$      & $29.31 $ &$ 29.25 $ & $\mathbf{29.81}$ \\
Mosaic      &MSSIM    & $0.9203$ & $0.8404$ &$0.8662$ &$0.9189$ &$0.9144$     & $0.9078$ &$ 0.7730$ & $\mathbf{0.9333}$ \\
            &PSNR     & $27.61$  & $26.93$  &$26.55$  &$27.66$  &$27.71$      & $27.34 $ &$ 26.51$  & $\mathbf{27.87}$ \\   
\bottomrule
\end{tabular}
}
\end{table}
\begin{table}[]
\centering
\vspace{1.0cm}
\caption{Comparison of ISNR and $\text{PSNR}_{\text{grad}}$ values for various approaches and proposed model. Clean image is degraded by additive Gaussian noise of mean $0.0$ and $\sigma=40$.}
\label{tab:ISNRtable}
        \scalebox{0.6}{
        \begin{tabular}{cccccccccc}
\toprule
Images & Measure   &TV\cite{rudin1992nonlinear} &NS\cite{nitzberg1992nonlinear} &Luo\cite{luo2006coupled} &RD\cite{guo2011reaction} &NLM\cite{buades2005non} & PV\cite{prasath2014system}  & SYS\cite{sun2016class} &Proposed\\
\midrule
Boat    &ISNR                            &$7.02$ &$7.01$  &$6.93$ &$7.78$ &$9.90$  &9.99  &10.19 & $\mathbf{10.80}$ \\
        &$\text{PSNR}_{\text{grad}}$     &$27.04$ &$27.01$ &$26.86$ &$28.07$ &$30.84$ &30.90 &31.28 & $\mathbf{31.64}$ \\
Lake    &ISNR                            &$7.05$ &$7.12$  &$7.00$ &$7.71$ &$9.12$  &9.19  &9.43  & $\mathbf{10.00}$ \\
        &$\text{PSNR}_{\text{grad}}$     &$27.23$ &$27.30$ &$27.12$ &$28.05$ &$30.11$ &30.19 &30.21 & $\mathbf{30.62}$ \\
Livingroom  &ISNR                        &$7.51$ &$7.02$  &$7.22$ &$8.16$ &$9.03$  &9.20  &9.36  & $\mathbf{9.77}$ \\
        &$\text{PSNR}_{\text{grad}}$     &$27.68$ &$27.33$ &$27.51$ &$28.62$ &$29.84$ &30.04 &30.25 & $\mathbf{30.62}$ \\
Mandril     &ISNR                        &$7.12$ &$7.19$  &$7.31$ &$7.58$ &$7.77$  &8.51  &8.66  & $\mathbf{8.86}$ \\
        &$\text{PSNR}_{\text{grad}}$     &$26.59$ &$26.70$ &$26.83$ &$27.44$ &$27.91$ &28.48 &28.58 & $\mathbf{28.88}$ \\
Pirate      &ISNR                        &$7.92$ &$7.83$  &$7.64$ &$8.66$ &$9.96$  &10.14  &9.99 & $\mathbf{10.66}$ \\
       &$\text{PSNR}_{\text{grad}}$      &$28.11$ &$27.57$ &$27.41$ &$29.14$ &$30.59$ &30.62 &30.60 & $\mathbf{31.26}$ \\
Walkbridge  &ISNR                        &$6.59$ &$6.88$  &$6.62$ &$7.07$ &$7.30$  &7.55  &7.59  & $\mathbf{7.78}$ \\
           &$\text{PSNR}_{\text{grad}}$  &$26.55$ &$26.65$ &$26.59$ &$27.14$ &$27.68$ &27.71 &27.73 & $\mathbf{27.98}$ \\
Woman      &ISNR                         &$9.11$ &$10.95$ &$9.29$ &$9.97$ &$13.11$ &13.94 &13.87 & $\mathbf{15.01}$ \\
           &$\text{PSNR}_{\text{grad}}$  &$30.05$ &$32.36$ &$30.72$ &$31.46$ &$34.77$ &37.37 &37.30 & $\mathbf{38.80}$ \\
Mosaic     &ISNR                         &$12.57$ &$11.64$ &$11.51$ &$12.60$ &$12.74$ & 12.17 &11.29 & $\mathbf{12.86}$ \\
           &$\text{PSNR}_{\text{grad}}$  &$37.34$ &$35.68$ &$35.48$ &$37.37$ &$37.52$ &36.01 &34.92 & $\mathbf{37.56}$ \\
\bottomrule
\end{tabular}
}
\end{table}

In addition to qualitative comparisons, the quantitative results, in terms of PSNR and MSSIM values for different test images as well as noise levels ($\sigma=20, 30, 40, 50$), are shown in tables [\ref{tab:table20}]-[\ref{tab:table50}]. To ease the comparison, the highest values of both MSSIM and PSNR measures are highlighted in each table. The highest values of PSNR and MSSIM for each image, clearly shows that the proposed model is better than all other models considered in the table. We have a significant observation that proposed model performs better for higher noise densities, as can be seen from tables[\ref{tab:table40}]-[\ref{tab:table50}]. For low noise densities, such as $\sigma=20,$ our results are still batter or very closed to the best values. Hence, analysis of the results presented in the tables, reveals that the performance of the present model increases with increasing noise density.
In table [\ref{tab:ISNRtable}], apart from the MSSIM and PSNR, we show the quantitative comparison with gradient PSNR and ISNR values of the proposed model and the alternative approaches. Finally, we note that the different quantitative measures (in terms of PSNR, MSSIM, gradient PSNR, and ISNR values) confirm the superiority of the proposed approach in comparison with other models.
Our comparison study shows that the proposed CPDE model preserves image structure efficiently in comparison to other earlier reported methods when signal-to-noise ratio is low. On the other hand, when signal-to-noise ratio is high, the present model restores image structure better than the other existing models. From the above discussion, it is confirmed that the proposed model is robust and more efficient than the available techniques considered here.
\section{Conclusion}
\label{sec:conclusion}
In this work, we propose a new CPDE based denoising framework. The proposed space-time regularization based image denoising approach highlights the choice of diffusion function and data fidelity term. For this purpose, we use the evolution equations to determine both the terms in the proposed model. It may be treated as the non-linear diffusion model which yield two separate PDEs to remove the noise with preservation of significant edges and fine structures. First, we establish the well-posedness of the proposed model using a time discretization method. Then, to solve the proposed model numerically, an implicit finite difference scheme along with advanced iterative solver has been used. We validate our model with different standard test images. Experimental studies confirm that the proposed model is more efficient than the existing models, in terms of quantitative measures and visual quality. Especially, the present model produces best results for high noise levels. We believe that present model has multidimensional applications which will ultimately benefit the society.

\section*{Competing interests}
The authors declare that they have no competing interests.
\thispagestyle{empty}

\section*{Acknowledgments}
Authors sincerely acknowledge Dr. Arnav Bhavsar, School of Computing and Electrical Engineering, Indian Institute of Technology Mandi, Himachal Pradesh, India for some valuable discussions and suggestions during the preparation of this manuscript.

\bibliographystyle{unsrt}  


\end{document}